\documentclass[11pt]{article}
\usepackage{graphicx} 

\usepackage{xcolor}

\usepackage{amssymb}
\usepackage{amsmath}
\usepackage{amsthm}

\usepackage{authblk}

\usepackage{hyperref}

\usepackage{subfigure}

\newtheorem{theorem}{Theorem}[section]
\newtheorem{lemma}[theorem]{Lemma}

\theoremstyle{plain} 

\theoremstyle{definition} 
\newtheorem{definition}[theorem]{Definition}

\title{Parallel Newton–Krylov Methods for 3D Free-Surface Flows under Strong Capillarity: A Lagrange Multiplier Approach}

\author{T. Benkley\thanks{Laboratory of Mechanical Metallurgy, EPFL, Route Cantonale, Lausanne, 1015, Vaud, Switzerland. Corresponding author \href{mailto:tyler.benkley@epfl.ch}{\nolinkurl{tyler.benkley@epfl.ch}} } \qquad S. Deparis\thanks{MATH, EPFL, Route Cantonale, Lausanne, 1015, Vaud, Switzerland. \href{mailto:simone.deparis@epfl.ch}{\nolinkurl{simone.deparis@epfl.ch}}}
}

\date{\today}

\newcommand{\Ca}{\textit{Ca}}
\newcommand{\Rey}{\textit{Re}}
\renewcommand{\u}{\boldsymbol{u}}
\renewcommand{\v}{\boldsymbol{v}}
\newcommand{\bnabla}{\boldsymbol{\nabla}}
\newcommand{\x}{\boldsymbol{x}}
\newcommand{\X}{\boldsymbol{X}}
\newcommand{\n}{\boldsymbol{n}}
\renewcommand{\d}{\boldsymbol{d}}

\newcommand{\F}{\boldsymbol{F}}

\renewcommand{\t}{\boldsymbol{t}}

\newcommand{\s}{\boldsymbol{s}}

\newcommand{\keywords}[1]{%
    \par\medskip
    \noindent\textbf{Keywords:} #1
}

\usepackage[margin=1in]{geometry}

\usepackage[
backend=bibtex,
sorting=none
]{biblatex}

\begin{document}

\maketitle

\begin{abstract}

An arbitrary Lagrangian-Eulerian method based on an alternative weak formulation, inexact Newton scheme and physics-based preconditioner is proposed for the solution of 3D steady laminar free-surface flows under high surface tension. Specifically, the Young-Laplace equation is imposed as the distinguishing condition for the free surface in lieu of the classical impermeability condition, relocating surface-tension forces from off-diagonal entries to a diagonal block of the linear system. This facilitates the approximation of surface-tension forces within preconditioners. Impermeability is instead imposed, via the Lagrange-multiplier method, as a constraint on the velocity field.

The physics-based preconditioner relies on a block LU factorization of the Jacobian in the high surface-tension limit. Diagonal blocks of the preconditioner are inverted either by multi-grid algorithms or by spectrally equivalent operators. In particular, the impermeability Schur complement is shown to be spectrally equivalent to a fractional Laplace operator which is approximated via a rational function.

Numerical simulations of the die-swell problem and of a 3D coating flow demonstrate robustness of the preconditioner, opening the door to large-scale simulations in the surface-tension-dominated regime with the arbitrary Lagrangian-Eulerian approach. Specifically, results show Krylov-space iteration counts between 70 and 150 over four decades of the Capillary number and suggest mild sensitivity to the mesh size. Parallelized computations on a cluster enable simulations to be scaled up to 68 million degrees of freedom in certain cases, while achieving a total runtime on the order of tens of minutes. Simulations furthermore exhibit convergence of the inexact Newton scheme in at most six iterations in the cases investigated here. Die-swell radius predictions are validated against earlier simulations and experiments.


\end{abstract}

\section{Introduction} \label{Sec:Introduction}

Surface-tension effects often defy intuition, leading to the formation of nontrivial flow patterns and complex liquid shapes which may sometimes only be explained by numerical simulation. At small scales, surface-tension forces take on a dominant role and can completely change the complexion of fluid flows, giving rise to unique phenomena such as wicking, wine tears over a glass surface or adhesion mechanisms for certain insects \cite{walker_adhesion_1993}. A scientific understanding of capillarity's effects has proven to be crucial, for instance, in the fine-tuning of micro-fluidic devices \cite{stone_engineering_2004}, in limiting smearing in ink-jet printing \cite{wijshoff_drop_2018} or in manufacturing methods such as welding \cite{kou_welding_2003}. 

Numerical algorithms for flows under surface tension must be capable of accounting for a dependence of the fluid flow upon an \textit{unknown} domain, which itself is a product of the flow, leading to a free-boundary problem. In such problems, one must solve for both the fluid-flow unknowns and the domain simultaneously, something that only a handful of specialized numerical methods are equipped to do. Many of these specialized techniques fall into the two main categories of (i) Eulerian and (ii) arbitrary Lagrangian-Eulerian. 

Eulerian methods discretize the numerical domain into a fixed grid of cells that the free interface is permitted to cut through \cite{scardovelli_direct_1999}. Rather than enforcing surface-tension forces over the free interface, they model capillarity as a volumetric force that is distributed over a neighborhood of the free interface \cite{brackbill_continuum_1992}. The thickness of this neighborhood is a user-prescribed quantity, which sets the scale for the approximation error. Naturally, the minimum thickness of this interface is limited from below by the resolution of the grid. A few acclaimed Eulerian approaches are the Volume-of-Fluid method \cite{hirt_volume_1981}, the Level-set method \cite{osher_fronts_1988}, the Phase-Field technique \cite{cahn_free_1958,kim_phase-field_2012}, and the Immersed-Boundary method \cite{peskin_immersed_2002}. While these methods are well equipped to handle topological changes such as breakup or coalescence, in high-capillarity flows, their surface-tension approximation may result in spurious velocity currents, even in simple situations such as that of the static Laplace droplet. The root cause has been identified as an inability to devise a volumetric approximation of surface-tension forces that is simultaneously "well-balanced" and consistent (see \cite{popinet_numerical_2018} for a thorough review of the topic). 

On the other hand, Arbitrary Lagrangian-Eulerian (ALE) methods \cite{stein_arbitrary_2004}, the main focus here, employ front-tracking. To be more specific, instead of implicitly tracking the free interface through a volume fraction or level set function, an explicit representation of the free surface is kept at all times. Thus, ALE approaches need not smear surface tension forces over a finite volume of high resolution, but rather can enforce them only over the free surface of interest, thereby side-stepping errors native to the approximation and the issue of spurious currents. This effectively makes them more accurate than Eulerian counterparts at dealing with problems devoid of topological change and of large mesh deformation, which are the primary focus here. 

Having been applied to myriad free-boundary problems \cite{saito_study_1981,kistler_coating_1984,ettouney_finite-element_1983,ruschak_method_1980,cuvelier_thermocapillary_1986,cuvelier_numerical_1990,sackinger_newtonraphson_1996,cairncross_finite_2000,baer_finite_2000,sprittles_finite_2012,anthony_sharp_2023}, Newton's algorithm has emerged as a robust and popular solution method for ALE problems. It owes its robustness to a simultaneous consideration, in Jacobian computations, of fluid flow and domain sensitivities. It suffers, however, from the drawback that, at each nonlinear iteration, a large monolithic linear system encompassing both the fluid-flow unknowns and the domain must be solved. Moreover, these linear systems have been reported to be poorly conditioned, especially in surface-tension-dominated problems, something that has inhibited the use of Krylov-space methods \cite{schunk_iterative_2002,osti_1372203}. To overcome these issues, research efforts have been invested in developing preconditioners to remedy the poor conditioning of ensuing linear systems and, in parallelized algorithms, to cope with a potentially larger problem size.

A group at Sandia National Laboratories attempted monolithic incomplete LU factorization (ILU) schemes \cite{cairncross_finite_2000,baer_finite_2000}, as well as coupling thresholded ILU with global reordering \cite{schunk_iterative_2002}. Unfortunately, these preconditioning strategies relied on a stabilized weak formulation to ensure convergence of the iterative solver, thereby trading numerical accuracy for improved conditioning \cite{osti_1372203}. Even with stabilization, however, ILU preconditioning has been reported to fail at addressing conditioning issues in certain cases \cite{osti_1372203}. Furthermore, the loss of accuracy associated with stabilization is particularly pronounced in the surface-tension-dominated regime \cite{osti_1372203}. These preconditioners were also found to exhibit limited parallel scalability when coupled with non-overlapping domain decomposition \cite{Baer1999ParallelSO}.

Another approach, in \cite{zhou_assessment_2001}, was diagonal preconditioning with a restarted Generalized Minimum Residual method in the context of massively parallel three-dimensional flows; however, this was lamented as "fragile" by the authors.

In the context of viscoelastic flows, Castillo \textit{et al.} proposed a parallel monolithic sparse approximate inverse preconditioner \cite{castillo_parallel_2009}. This method outperformed ILU preconditioning on a viscoelastic flow example, also exhibiting parallel scalability. Its algorithmic scalability, however, was not assessed in the scope of \cite{castillo_parallel_2009} and, furthermore, large Krylov subspaces were required to reach convergence.

In \cite{benkley_generalization_2026}, a preconditioner based on a block LU factorization was implemented to solve a fluid-flow problem under surface tension. Specifically, the mixed finite element system was decomposed into variable-specific blocks and factorized into an upper block-triangular matrix, with each diagonal block preconditioned according to the structure of its associated operator. Numerical results were presented for serial 3D simulations and verified on a problem where the Krylov subspace sizes required to reach convergence were relatively small. The performance of the preconditioner, however, was explored predominantly in cases where surface-tension values were moderate, with results suggesting a decrease in effectiveness at larger values. It also remains to be investigated whether this technique can scale in parallel.

Thus, the downsides of system size and ill-conditioning have unfortunately hindered the development of scalable ALE algorithms, particularly in surface-tension-dominated problems, and, as a result, their scope of applicability has been limited to 2D or small 3D problems.

In this article, an ALE method is proposed for the parallelized solution of the 3D steady Navier-Stokes equations in the surface-tension-dominated regime for low Reynolds numbers. Specifically, in Section \ref{Sec:NavierStokes}, after defining the classical weak formulation of the Navier-Stokes equations with a free capillary surface, an alternative weak form of the problem is given. This formulation relies on the imposition of the Young-Laplace equation as the distinguishing condition for the free-surface position instead of the impermeability condition, which is imposed via the Lagrange-multiplier method. This formulation is favorable for the high surface-tension regime because capillary forces emerge through a diagonal block, as opposed to existing approaches \cite{benkley_generalization_2026}. In Section \ref{Sec:NumericalMethod}, an inexact Newton method is presented for the free-boundary-problem solution. Via an apt choice of the displacement closure equations, within preconditioners, displacement degrees of freedom over the free surface are decoupled from those in the bulk. Then, after a brief summary of the Finite Element Method (FEM) used here, the linear systems ensuing from a discretization by the FEM are provided in block formulation. In Section \ref{Sec:Preconditioning}, a preconditioner based on a block LU factorization of the Jacobian in the high surface-tension limit is derived. In addition, the Schur complement of the impermeability multiplier is proven to be spectrally equivalent to an operator inducing a fractional Sobolev norm. Block diagonals are then substituted by spectrally equivalent operators or a multi-grid cycle. In Section \ref{Sec:NumericalSimulations}, the preconditioner and the ALE formulation are implemented to solve two 3D numerical examples, namely (i) the die-swell problem and (ii) a coating-flow problem. In the case of the die swell, the present method is validated against existing experiments \cite{goren_shape_1966} and computations by an earlier ALE formulation \cite{benkley_generalization_2026}. Coating-flow results across four orders of magnitude of the Capillary number exhibit widely different free surfaces, qualitatively consistent with the physics of wetting. In both the die swell and coating flow, the performance of the preconditioner is assessed through the number of GMRes iterations. Nonlinear convergence of the inexact Newton scheme is furthermore evaluated. Clock times are reported for the die-swell for two different levels of refinement with the finest simulations reaching over 68 million degrees of freedom. Typical simulation meshes are shown in Appendices \ref{Sec:AppendixA} and \ref{Sec:AppendixB}. Results are then discussed and rationalized in Section \ref{Sec:Discussion}.

\section{The Navier-Stokes Equations with a Free Capillary Surface} \label{Sec:NavierStokes}

In this section, after first introducing the requisite function spaces and the mathematical model of steady incompressible fluid flows with a free capillary surface, an alternative weak form of surface-tension-dominated flows is proposed.

The problem of interest here is that of steady, incompressible fluid flows where the liquid forms a free interface with a quiescent, gas atmosphere. Small-scale flows are considered, implying that surface-tension forces are dominant. Such problems are modeled by the incompressible Navier-Stokes equations where the free interface is governed by the Young-Laplace equation, a shear-free condition and an impermeability requirement.

Let the liquid domain be $\Omega\subset \mathbb{R}^3$, where $\partial \Omega$ is formed by two disjoint surfaces, $\Sigma_r$, the rigid surface, and $\Sigma_f$, the free surface. Surface $\Sigma_r$ is furthermore defined as the union of two disjoint surfaces, $\Sigma_{N,r}$ and $\Sigma_{D,r}$. Let $\Gamma_c = \partial\Sigma_f\cap \partial \Sigma_{N,r}$ be the contact line, that is the three-phase line forming an intersection between liquid, solid and gas. 

Moreover, let us define the following Sobolev spaces and respective inner products
\[ L^2(\Omega) = \{q:\Omega\rightarrow \mathbb{R} \mid \|q\|_{L^2(\Omega)} < \infty \},\]

\[\langle p,q \rangle_{L^2(\Omega)} = \int_\Omega p\,qdV,\textrm{ for $p,q\in L^2(\Omega)$},\]

\[ [H^1(\Omega)]^d = \{ \v:\Omega\rightarrow \mathbb{R}^d \mid \|\v\|_{[H^1(\Omega)]^d} < \infty \},\]

\[\langle \u,\v \rangle_{[H^1(\Omega)]^d} = \int_\Omega \u\cdot\v dV + \int_\Omega \bnabla\u:\bnabla\v dV, \textrm{ for $\u,\v\in [H^1(\Omega)]^d$},\]

\[H^{1/2}(\Sigma) = \{ f:\Sigma \rightarrow \mathbb{R} \mid \Sigma\subseteq \partial \Omega,\,\|f\|_{H^{1/2}(\Sigma)}<\infty \},\]

\[\langle f,g \rangle_{H^{1/2}(\Sigma)} = \int_\Sigma (1-\Delta_S)^{1/4}f\, (1-\Delta_S)^{1/4}g  dS,\textrm{ for $f,g\in H^{1/2}(\Sigma)$},\]
where $\Delta_S$ is the Laplace-Beltrami operator. The dual space of a function space $X$ is denoted by $X'$. The inner product of $(H^{1/2}(\Sigma))'$ is defined as
\[\langle \chi,\lambda \rangle_{(H^{1/2}(\Sigma))'} = \int_\Sigma (1-\Delta_S)^{-1/4}\chi\, (1-\Delta_S)^{-1/4}\lambda  dS,\textrm{ for $\chi,\lambda \in (H^{1/2}(\Sigma))'$}. \]
The definitions of $H^{1/2}(\Sigma)$'s inner product and that of its dual space follow the theory of fractional power spaces in \cite{lions_non-homogeneous_1972}. 

\begin{figure}
    \centering
    \subfigure[]{\includegraphics[width=0.49\linewidth]{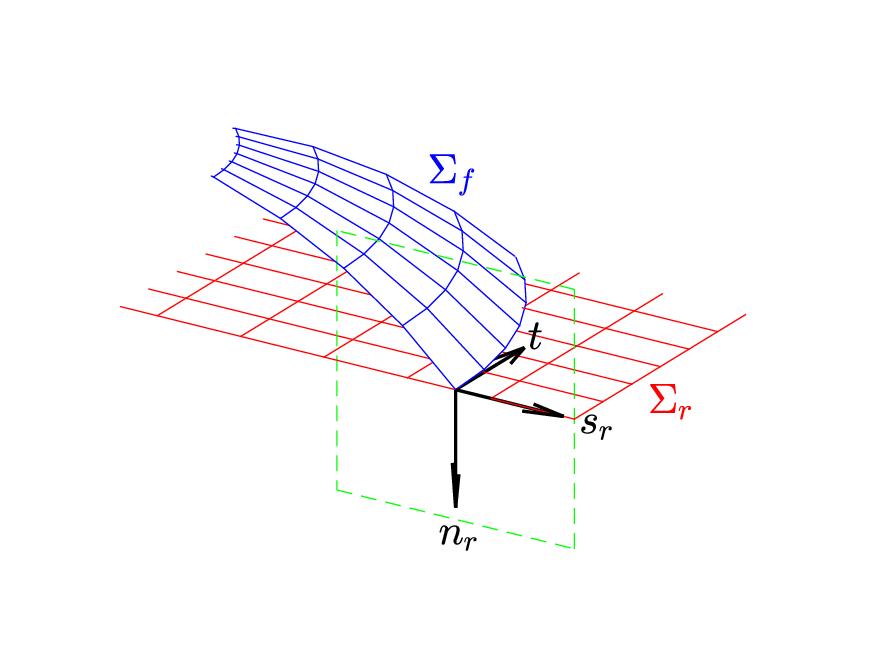}}
    \subfigure[]{\includegraphics[width=0.49\linewidth]{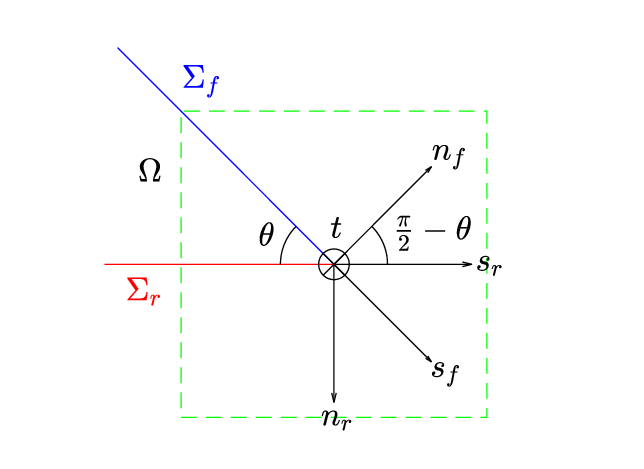}}
    \caption{Orthonormal bases at the contact line; Adapted from Benkley \textit{et al.} \cite{benkley_generalization_2026}, under the CC BY 4.0 license.}
    \label{fig:ContactLines}
\end{figure}

Let the function spaces $V$ and $Q$ respectively be $[H^1(\Omega)]^3$ and $L^2(\Omega)$. Let $\u\in V$, $p\in Q$ and $p_a\in \mathbb{R}$ be the velocity field, pressure field and atmospheric pressure constant respectively, and let $\tau$ be defined as the Newtonian stress tensor, \textit{i.e.} $\tau = -p\boldsymbol{I} +\mu ( \nabla \u + (\nabla \u)^T)$, $\boldsymbol{I}$ being the identity tensor. The space $S$ is assumed to be a suitable space for the domain $\Omega$. In this case, the physical problem is modeled by

\begin{equation} \label{Eq:StrongForm}
    \begin{split}
    &\text{Find $(\u,p,\Omega,p_a)\in V\times Q \times S \times \mathbb{R} $ such that} \\
    &\begin{cases}
        \rho (\u\cdot \bnabla)\u = -\bnabla p + \mu \nabla^2\u, \hfill \text{ in $\Omega$} \\
        \bnabla \cdot \u = 0, \hfill\text{ in $\Omega$} \\
        \u\cdot\n = 0, \hfill\text{ over $\Sigma_f$} \\
        \tau \n = -\gamma(\bnabla_S \cdot \n)\n - p_a\n, \hfill\text{ over $\Sigma_f$} \\
        \n_f\cdot \n_r = -\cos(\theta), \hfill\text{ along $\Gamma_c$}\\
        \tau\n = \boldsymbol{g_N},\hfill\text{ over $\Sigma_{N,r}$} \\
        \u = \boldsymbol{u_D},\hfill\text{ over $\Sigma_{D,r}$} \\
        \text{+ global constraint}
    \end{cases}
    \end{split}
\end{equation}
where $\rho$, $\mu$ and $\gamma$ are the density, dynamic viscosity and surface tension, respectively, all of the liquid. The parameter $\theta$ denotes the equilibrium contact angle, a physical property for a given solid, liquid and gas. Moreover, $\n$ denotes the unit normal vector of the surface over which it is evaluated. Along the curves at the intersection of surfaces, where the notion of normal vector is ambiguous, the subscripts $f$ and $r$ are employed to distinguish between the normal vectors of $\Sigma_f$ and $\Sigma_r$ respectively, \textit{i.e.} $\n_f = \n|_{\Sigma_f}$ and $\n_r = \n|_{\Sigma_r}$ (see Figure \ref{fig:ContactLines}). The symbol $\bnabla_S$ is defined as the surface gradient, \textit{i.e.} $\bnabla_S = (\boldsymbol{I} - \n\n^T)\bnabla$. The functions $\boldsymbol{g_N}\in (H^{1/2}(\Sigma_{N,r}))'$ and $\u_D\in H^{1/2}(\Sigma_{D,r})$ are prescribed, respectively, for the traction vector over $\Sigma_{N,r}$ and the velocity field over $\Sigma_{D,r}$.

The global constraint, associated with the atmospheric pressure constant, depends on the problem of interest. In enclosed flows, one may use a volume constraint while in inflow-outflow problems, a constraint on the outflow-surface radius or area may be preferred.

The following variables are non-dimensionalized as follows:
\begin{align}
\x^* := \frac{\x}{L}, && \u^* := \frac{\u}{U}, && p^* := \frac{Lp}{\mu U}, && p_a^* := \frac{L p_a}{\gamma},
\end{align}
where $\x\in\mathbb{R}^3$ denotes the position in space. $L$ and $U$ are problem-dependent, characteristic scales chosen for the length and the velocity respectively. Non-dimensionalizing Equation (\ref{Eq:StrongForm}) with these characteristic scales gives rise to the dimensionless numbers,
\begin{align}
\Rey := \frac{\rho U L}{\mu}, && \Ca := \frac{\mu U}{\gamma},
\end{align}
respectively, the Reynolds and Capillary numbers. Hereafter, the superscripts $*$ are omitted and all quantities are made to be dimensionless.

The standard approach when solving the Navier-Stokes equations with a free surface is to implement the impermeability condition as the \textit{distinguishing condition} for the free surface (as defined in \cite{cairncross_finite_2000}), \textit{i.e.} the equation associated with the test function of the free-surface unknown.  The Young-Laplace equation is then imposed as a natural boundary condition in the momentum equation. Integration by parts is employed using the divergence theorem for curved 3D surfaces \cite{weatherburn_differential_1927} to produce the subsequent dimensionless weak form

\begin{equation} \label{Eq:ClassicalWeakForm}
    \begin{split}
    &\text{Find $(\u,p,\Omega,p_a)\in V_D\times Q \times S \times \mathbb{R} $ such that} \\
    &\begin{cases}
        \int_\Omega \{\,\Rey\,\v\cdot (\u\cdot \bnabla )\u  + [\bnabla\u+ (\bnabla\u)^T]:\bnabla\v  -p(\bnabla\cdot\v)\}dV \\
        \text{\textcolor{white}{aaaaa}}= - \frac{1}{\Ca}\int_{\Sigma_f}\{\bnabla_S\cdot\v + p_a(\v\cdot\n)\}dS \\
        \text{\textcolor{white}{aaaaa=/}}+  \frac{1}{\Ca}\int_{\Gamma_c} \cos{(\theta)}(\v\cdot\s_r) dl \\
        \text{\textcolor{white}{aaaaa=/}}+  \int_{\Sigma_{N,r}} (\v\cdot \boldsymbol{g_N}) dS,\,\hfill \forall \v\in V_0 \\
        \int_\Omega -q(\bnabla \cdot \u)dV = 0,\,\hfill\forall q\in Q \\
        \int_{\Sigma_f} \chi(\u\cdot \n)dS = 0,\,\hfill\forall \chi \in L \\
        \text{+ global constraint}
    \end{cases}
    \end{split}
\end{equation}
where $V_D =\{ \v\in V \mid \v|_{\Sigma_{D,r}} = \u_D \}$, $V_0 = \{ \v\in V \mid \v|_{\Sigma_{D,r}} = 0 \}$ and $L= (H^{1/2}(\Sigma_f))'$. The variables, $\v\in V_0$, $q\in Q$ and $\chi\in L$ are test functions. Furthermore, as shown in Figure \ref{fig:ContactLines}, $\s_r$ is the unique unit vector that is both normal to $\Gamma_c$ while being tangent to $\Sigma_r$ and pointing outward. 
The function space $V_D$ is equipped with the energy inner product \[ \langle \u,\v\rangle_{V_D} = \frac{1}{2}\int_\Omega [\bnabla\u + (\bnabla \u)^T]:[\bnabla \v+(\bnabla\v)^T] dV,\, \forall \u,\v\in V_D, \] which induces an inner product equivalent to the standard $[H^1(\Omega)]^3$ one for the flows studied here. The norm used here is induced by the inner product.

If capillarity is imposed weakly as it is in Equation (\ref{Eq:ClassicalWeakForm}), this causes Young-Laplace terms to be located in off-diagonal blocks of the Jacobian matrix. In preconditioners based on a block LU factorization such as those studied in \cite{benkley_generalization_2026}, capillarity thus manifests itself through a dense Schur complement which, ideally, must be approximated by a sparse operator. Devising an accurate sparse approximation can be challenging. Preconditioner performance then hinges on the quality of this approximation.

The scheme proposed next, instead, bypasses the need to approximate this Schur complement by enforcing the Young-Laplace equation as the distinguishing condition for the free surface. To this end, the impermeability condition is employed instead as a boundary condition on the velocity field. Different methods exist for imposing such a boundary condition over a curved boundary, such as Nitsche's method \cite{Nitsche1971berEV,juntunen_nitsches_2008} or the Lagrange multiplier technique \cite{babuska_finite_1973,Verfürth1986/87}, the latter being used here. Let $\lambda\in L$ be the Lagrange multiplier and $\chi\in L$ its associated test function. The ensuing weak form may be found below.

\begin{equation} \label{Eq:WeakForm}
    \begin{split}
    &\text{Find $(\u,p,\lambda,\Omega,p_a)\in V_D\times Q \times L_0 \times S \times \mathbb{R} $ such that} \\
    &\begin{cases}
        \int_\Omega \{\,\Rey\,\v\cdot (\u\cdot \bnabla )\u  + [\bnabla\u+ (\bnabla\u)^T]:\bnabla\v  -p(\bnabla\cdot\v)\}dV - \int_{\Sigma_f} \lambda (\v\cdot\n) dS \\
        \text{\textcolor{white}{aaaaa}}=  \int_{\Sigma_{N,r}} (\v\cdot \boldsymbol{g_N}) dS,\,\hfill\forall \v\in V_0 \\
        \int_\Omega -q(\bnabla \cdot \u)dV = 0,\,\hfill\forall q\in Q \\
        \int_{\Sigma_f} -\chi(\u\cdot \n)dS = 0,\,\hfill\forall \chi\in L_0 \\
        \int_{\Sigma_f} [ (\boldsymbol{\xi}\cdot\n)( \Ca \,\lambda + p_a) + \bnabla_S\cdot\boldsymbol{\xi}  ] dS - \cos{(\theta)} \int_{\Gamma_c} (\boldsymbol{\xi}\cdot \s_r)= 0,\,\hfill\forall \boldsymbol{\xi}\in M \\
        \text{+ global constraint}
    \end{cases}
    \end{split}
\end{equation}
where $M$ is a suitable function space defined in the next section.

Let us define the function \(\zeta := \v\cdot\n \) over \(\Sigma_f\) for a given \(\v\in V_D\). Since \(V_D\subset[H^1(\Omega)]^3\), \(\zeta\) belongs to \(H^{1/2}(\Sigma_f)\). Under suitable regularity assumptions on \(\Sigma_f\), the boundary conditions of \(\v\) over \(\Sigma_{D,r}\) furthermore permit the zero extension of \(\zeta\),
\[
\widetilde\zeta=
\begin{cases}
\zeta, &\textrm{over }\Sigma_f,\\
0, &\textrm{over }\Sigma_{D,r}.
\end{cases}
\]
Denoting the intersection of $\partial \Sigma_f$ and $\partial \Sigma_{D,r}$ by \(\Gamma_D:=\partial\Sigma_f\cap\partial\Sigma_{D,r}\), following \cite{lions_non-homogeneous_1972}, the Lions–Magenes space is defined as
\[
H_{00,\Gamma_D}^{1/2}(\Sigma_f)
:= \{ \zeta\in H^{1/2}(\Sigma_f) \mid \widetilde\zeta\in H^{1/2}(\Sigma_f\cup\Sigma_{D,r})\}.
\]
As \(\lambda\) acts on \(\v\cdot\n\) through a duality pairing over \(\Sigma_f\), \(L_0\) is defined as the dual space of \(H_{00,\Gamma_D}^{1/2}(\Sigma_f)\). Under suitable regularity assumptions, the inner product, 
\[
\langle \chi, \lambda \rangle_{L_0} = \int_{\Sigma_f} (-\Delta_{S,0})^{-1/4}\chi (-\Delta_{S,0})^{-1/4}\lambda dS, \textrm{ for } \chi,\lambda\in L_0,
\]
induces an equivalent norm in $(H^{1/2}_{00,\Gamma_D}(\Sigma_f))'$ \cite{bonito_numerical_2013,chandler-wilde_interpolation_2015}. Here, \(\Delta_{S,0}\) denotes the Laplace–Beltrami operator satisfying homogeneous Dirichlet boundary conditions along \(\Gamma_D\) and homogeneous Neumann boundary conditions on the complement of \(\partial\Sigma_f\). In the problems studied here, \(\Gamma_D\) has positive boundary measure and \(\Sigma_f\) is connected, ensuring that \(-\Delta_{S,0}\) is strictly positive and therefore invertible.


Variable $\lambda$ has a physical meaning, as Verfürth showed that $\lambda = \n\cdot(\tau\n)$ over $\Sigma_f$ \cite{Verfürth1986/87}. This makes it possible to replace the normal stress term not only in the momentum equation, but also in the Young-Laplace law. Substitution is in fact preferred mathematically since the trace of $\tau\in L^2(\Omega)$ is not defined.

As mentioned in the Introduction, Equation (\ref{Eq:WeakForm}) is a free-boundary problem, implying that the domain, $\Omega$, is part of the solution. Standard numerical methods for partial differential equations (PDEs) must be adapted, as they are predicated on the domain being known. Specialized techniques are therefore required to effect a solution in this case. In the next section, the numerical method for free-boundary problems used in this article is explained.

\section{Numerical Method} \label{Sec:NumericalMethod}

\subsection{Newton's Method for Free-Boundary Problems}

When solving free-boundary problems in their weak form, one is faced with the task of integrating the various equations over the domain, something that is infeasible by a numerical algorithm since $\Omega$ is unknown. As a result, the standard approach is to parametrize the weak form integrals by a reference domain $\Omega_0$. Specifically, one assumes the existence of a displacement field $\d : \Omega_0\rightarrow \mathbb{R}^3$ that deforms $\Omega_0$ into $\Omega$, \textit{i.e.}
\begin{equation}
    \Omega = \{ \X + \d(\X) \mid \X \in \Omega_0 \}.
\end{equation}
This permits one to recast weak form terms into integrals over $\Omega_0$ (or subsets) where $\d$ emerges through the nonlinear parametrized quantities in Equation (\ref{Eq:Parametrization}), thereby transforming the task of finding $\Omega$ into one of finding $\d$. 

\begin{equation} \label{Eq:Parametrization}
    \begin{cases}
        f(\x) = f(\X+\d(\X)) \\
        \bnabla_x = \F^{-T}\bnabla_X \\
        \n = \frac{\F^{-T}\n_0}{|\F^{-T}\n_0|} \\
        \t = \frac{\F\t_0}{|\F\t_0|} \\
        \int_{\Omega}.\,dV = \int_{\Omega_0}.\,\textrm{det}(\F) dV \\
        \int_{\Sigma}.\,dS = \int_{\Sigma_0}.\,|\F^{-T}\n_0|\textrm{det}(\F) dS \\
        \int_{\Gamma_c}.\,dl = \int_{\Gamma_{c,0}}.\,|\F \t_0| dl \\
    \end{cases}
\end{equation}
where $\F = \boldsymbol{I} + \bnabla\d$ denotes the deformation gradient tensor. Displacement therefore emerges in volume integrals despite the fact that no physical law governs it within the bulk of $\Omega_0$, leaving it unconstrained. Tangent components of $\d$ over $\Sigma_{f,0}$ are also unconstrained. Additional equations must therefore be supplied to close the mathematical problem. Here, for closure, a Laplace problem is imposed on $\d$ within the bulk of $\Omega_0$, and a projected Laplace-Beltrami problem is enforced on $\d$ over $\Sigma_{f,0}$.
\begin{equation} \label{Eq:Closure}
    \begin{cases}
        - \Delta \d = 0,\textrm{ in }\Omega_0 \\
        - (\boldsymbol{I} - \frac{\F_0^{-T}\n_0\n_0^T\F_0^{-1}}{|\F_0^{-T}\n_0|^2})\Delta_{S}\d = 0,\textrm{ over }\Sigma_{f,0} ,
    \end{cases}
\end{equation}
where $\F_0$ denotes the deformation gradient tensor with argument $\d_0$, the current guess of the displacement field. The projection of the Laplace-Beltrami operator onto tangent directions of the deformed free surface ensures that, at convergence, its weak form is not evaluated over normal components of the test function that are used to test the Young-Laplace equation. This prevents a non-physical coupling between the Young-Laplace equation and Laplace-Beltrami equation from occurring. The boundary conditions imposed for the closure equations are problem-dependent and mentioned in Section \ref{Sec:NumericalSimulations}.

This particular choice of closure equations ensures that displacement degrees of freedom over the free surface are decoupled from interior degrees of freedom within $\d$'s block diagonal of the preconditioner.

Let us note that the Laplace-Beltrami equation in (\ref{Eq:Closure}) is nonlinear without being the Jacobian of a nonlinear problem, making this scheme an inexact Newton method. Projecting the Laplace-Beltrami equation over the deformed tangents from the previous iteration enables one to bypass the need to compute second-order derivatives in space of $\d$, which arise when taking the exact Jacobian of
\[ -(\boldsymbol{I} - \frac{\boldsymbol{F}^{-T}\n_0\n_0^T\boldsymbol{F}^{-1}}{|\boldsymbol{F}^{-T}\n_0|^2})\Delta_S\d = 0, \textrm{ over $\Sigma_{f,0}$}.\]

Considering that both the Young-Laplace equation and the Laplace-Beltrami problem are second-order PDEs of $\d$ over $\Sigma_{f,0}$, after integration by parts, the trace of $\d$ must satisfy the condition $\d|_{\Sigma_{f,0}} \in [H^1(\Sigma_{f,0})]^d$, which is necessary for well-posedness. The function space $M$ is thus defined as 
\[M = \{ \d\in [H^1(\Omega_0)]^3 \mid \d|_{\Sigma_{f,0}} \in [H^1(\Sigma_{f,0})]^3 \},\]
and equipped with the norm
\[\|.\|^2_M = \|.\|^2_{H^1(\Omega_0)} + \|.\|_{H^1(\Sigma_{f,0})}^2.\]
Furthermore, the space 
\[
M_\textrm{bulk} = \{ \d\in M \mid \d|_{\Sigma_{f,0}} = 0\}
\]
is defined.

The discretized formulation then becomes amenable to solution by a nonlinear algorithm. In this article, Newton's method is executed by computing and inverting the Jacobian of the continuous weak form at each nonlinear iteration. Contrary to the boundary-fitted approach which updates the domain guess $\Omega_0$ between nonlinear iterations, the entire computation is effected in a fixed reference domain, the reason being the following.

When confronted with a continuous set of equations, differentiation and discretization do not commute if the underlying domain, and in turn, the discrete function spaces, depend on the solution, something that occurs when the domain is updated between iterations. To avoid this potential source of inexactness of Newton's algorithm when using the Jacobian of the \textit{continuous} weak form, the numerical domain is kept constant, which renders the discrete function spaces independent of the solution. Meshes are then displaced upon convergence of the nonlinear iteration.

The nonlinear, parametrized weak form may be found below. 
\begin{equation} \label{Eq:ParametrizedWeakForm}
    \begin{split}
    &\text{Find $(\u,p,\lambda,\d,p_a)\in V_D\times Q \times L_0 \times M \times \mathbb{R} $ such that} \\
    &\begin{cases}
        \int_{\Omega_0} \{\,\Rey\,\v\cdot [\u\cdot (\F^{-T}\bnabla) ]\u  + [(\F^{-T}\bnabla)\u + ((\F^{-T}\bnabla)\u)^T ]:[(\F^{-T}\bnabla)\v]  \\
        \text{\textcolor{white}{aaaa}}- p[(\F^{-T}\bnabla)\cdot\v]\}\textrm{det}(\F)dV - \int_{\Sigma_{f,0}} \lambda [\v\cdot(\F^{-T}\n)]\textrm{det}(\F) dS \\
        \text{\textcolor{white}{aaaaaaaaaa}}=  \int_{\Sigma_{N,r,0}} (\v\cdot \boldsymbol{g_N}) |\F^{-T}\n|\textrm{det}(\F) dS,\,\hfill\forall \v\in V_0 \\
        \int_{\Omega_0} -q[(\F^{-T}\bnabla)\cdot\u]\textrm{det}(\F) dV = 0,\,\hfill\forall q\in Q \\
        \int_{\Sigma_{f,0}} -\chi[\u\cdot (\F^{-T}\n)]\textrm{det}(\F)dS = 0,\,\hfill\forall \chi\in L_0 \\
        \int_{\Omega_0} \bnabla\d:\bnabla\boldsymbol{\xi} dV = 0,\,\hfill\forall \boldsymbol{\xi}\in M_\textrm{bulk} \\
        \int_{\Sigma_{f,0}} \{ [\boldsymbol{\xi}\cdot(\F^{-T}\n)]( \Ca \,\lambda + p_a) +  [(I- \frac{(\F^{-T}\n)(\n^T\F^{-1})}{|\F^{-T}\n|^2})\F^{-T}\bnabla] \cdot\boldsymbol{\xi}  \}|\F^{-T}\n|\textrm{det}(\F) dS \\
        \text{\textcolor{white}{aaa}}- \cos{(\theta)} \int_{\Gamma_{c}} [\boldsymbol{\xi}\cdot (\F^{-T}\s_r)] \textrm{det}(\F) dl \\
        \text{\textcolor{white}{aaa}}+ \int_{\Sigma_{f,0}} \bnabla_S \d:\bnabla_S[(I - \frac{\F_0^{-T}\n_0\n_0^T\F_0^{-1}}{|\F_0^{-T}\n_0|^2})\boldsymbol{\xi}] dS = 0,\,\hfill\forall \boldsymbol{\xi}\in M \\
        \text{+ global constraint}
    \end{cases}
    \end{split}
\end{equation}
This equation system is discretized with the FEM, leading to a mixed FEM block formulation. The discretization is discussed next.

\subsection{Finite Element Method} \label{Sec:FEM}

The FEM, used to solve the numerical examples of Section \ref{Sec:NumericalSimulations}, is a numerical method designed to compute Galerkin approximations of PDE solutions. It discretizes the domain $\Omega_0$ of the PDE into a tesselation of non-overlapping elements which may be either tetrahedra (triangles in 2D) or hexahedra (quadrilaterals in 2D). Each element, of a characteristic user-prescribed mesh size $h$, is then provided with a finite number of nodes where the positioning and number of nodes depend on the type and degree of the finite element function space. 

The Lagrange finite elements, used here, are continuous, piecewise polynomials of degree $r\geq 1$ whose nodes are distributed uniformly along the edges of each element. A Lagrange finite element space of degree $r$, characteristic mesh size $h$ and domain $\Omega$ is denoted hereafter by $\chi_{r,h}(\Omega)$.

The space $\chi_{r,h}(\Omega)$ is then characterized by the basis $\{\phi_1,\phi_2,...,\phi_{N-1},\phi_N\}$ where each element of the basis $\phi_i$, defined over $\Omega$, is identically 1 at its respective node, and 0 at all others. Basis functions are therefore nonzero only within elements adjoined by their respective node, giving them a compact support. 

Given a basis of functions, an approximate solution $z_h \in \chi_{r,h}(\Omega)$ can then be expanded as a linear combination of basis elements.
\begin{equation}
    z_h =  \sum_{i=1}^{N} z_{h,i} \phi_i 
\end{equation}
Given a function space $Z$, its dual $Z'$, a bilinear form $a:Z\times Z\rightarrow \mathbb{R}$ and a linear functional $F\in Z'$, a generic PDE in weak form such as
\begin{equation} \label{Eq:GenericPDE}
\begin{split}
&\textrm{Find $z\in Z$ such that}\\
&a(z,y) = F(y),\,\forall y\in Z, 
\end{split}
\end{equation}
may be discretized by injecting the approximate solution $z_h\in \chi_{r,h}(\Omega)$ and testing the equation against each element of the basis of $\chi_{r,h}(\Omega)$.
\begin{equation}
\begin{split}
&\textrm{Find $z_h  \in \chi_{r,h}(\Omega)$ such that}\\
&\sum_{i=1}^N z_{h,i} a(\phi_i,\phi_j) = F(\phi_j),\,\forall j=1,...,N 
\end{split}
\end{equation}
From this discretization procedure ensues the linear system
\begin{equation} \label{Eq:GenericLinearSystem}
    A\vec{z_h} = \vec{f},
\end{equation}
where $A\in \mathbb{R}^{N \times N}$, $\vec{z_h}\in \mathbb{R}^N$, $\vec{f}\in \mathbb{R}^N$ and are defined by
\begin{equation}
\begin{split}
    &A_{ji} = a(\phi_i,\phi_j) \textrm{ for $i,j=1,...,N$} \\
    &\vec{z_h} = (z_{h,1},...,z_{h,N}) \\
    &f_j = F(\phi_j) \textrm{ for $j=1,...,N$}.
\end{split}
\end{equation}
Due to the fact that basis elements $\phi_i$ have a compact support, $A$ is a sparse matrix.

\subsubsection{Discretization}

In saddle-point PDEs such as the incompressible Navier-Stokes, the Ladyzhenskaya-Babuska-Brezzi (LBB) theorem, involving the inf-sup condition, must be satisfied in order to guarantee existence and uniqueness of a solution \cite{brezzi_existence_1974}. A classical theorem for two-variable saddle-point PDEs is:

\begin{theorem} \label{Thm:WellPosedness}
    \textbf{Ladyzhenskaya–Babuska–Brezzi} \\
    Let $X$ and $U$ be two Hilbert spaces, let $a:X\times X\rightarrow \mathbb{R}$ and $b:X\times U\rightarrow \mathbb{R}$ be two bilinear forms, and let $F \in X'$ and $G\in U'$ be two continuous functionals. The saddle-point problem
\[\begin{split}
    &\textrm{Find $(x,u)\in X\times U$ such that} \\
    &\begin{cases}
    a(x,y) + b(y,u) = F(y),\,\hfill\forall y\in X \\
    b(x,v) = G(v),\,\hfill\forall v\in U
    \end{cases}
    \end{split}
\]
is well-posed under the conditions that 
\begin{itemize}
    \item $a(.,.)$ and $b(.,.)$ are continuous bilinear forms, \textit{i.e.} there exist the continuity constants $\gamma_a$, $\gamma_b>0$ such that \[a(x,y) \leq \gamma_a \|x\|_X \|y\|_X \,\forall x,y\in X,\] 
    \[b(x,u) \leq \gamma_b \|x\|_X \|u\|_U \,\forall x\in X,\,\forall u\in U,\]
    \item $a(.,.)$ is coercive on the kernel of $b(.,.)$, \textit{i.e.} there exists a coercivity constant $\alpha>0$ such that \[a(x,x)\geq \alpha \|x\|_X^2,\,\forall x\in \textrm{ker}(b),\]
    \item $b(.,.)$ satisfies the inf-sup condition, \textit{i.e.} there exists an inf-sup constant $\beta>0$ such that \[ \underset{u \in U}{\textrm{inf}} \underset{x\in X}{\textrm{sup}} \frac{b(x,u)}{\|x\|_X\|u\|_U} \geq \beta, \]
\end{itemize}
and admits a unique solution $(x,u)\in X\times U$ satisfying the stability estimate
\[
\|x\|_X + \|u\|_U \leq C (\|F\|_{X'} + \|G \|_{U'})
\]
where the constant $C > 0$ depends only on $\gamma_a$, $\gamma_b$, $\alpha$, and $\beta$.
\end{theorem}
Discrete well-posedness furthermore requires that, for two finite-dimensional subspaces $X_h\subset X$ and $U_h\subset U$, the corresponding discrete stability constants $\gamma_{a,h}$, $\gamma_{b,h}$, $\alpha_h$ (defined on the discrete kernel $\textrm{ker}(b_h)$), and $\beta_h$ are bounded uniformly with respect to the mesh size $h$. Specifically, there must exist constants $\gamma_a^*, \gamma_b^*, \alpha^*, \beta^* > 0$, independent of $h$, such that $\gamma_{a,h} \leq \gamma_a^*$, $\gamma_{b,h} \leq \gamma_b^*$, $\alpha_h \geq \alpha^*$, and $\beta_h \geq \beta^*$ as $h \to 0$.

When effecting an approximate solution with the FEM, the inf-sup condition in Theorem \ref{Thm:WellPosedness} places restrictions on the stable combinations of finite element spaces. Although a proof of existence and uniqueness of a numerical solution of Equation (\ref{Eq:ParametrizedWeakForm}) falls outside of the scope of this article, consideration is given to the inf-sup condition when choosing the finite element spaces.

After discretizing the domain into a mesh of quadratic tetrahedra characterized by a mesh size $h$, finite element spaces for Problem (\ref{Eq:ParametrizedWeakForm}) are chosen to be $\u_h \in [\chi_{2,h}(\Omega_0)]^3,\,p_h \in \chi_{1,h}(\Omega_0),\,\lambda_h \in \chi_{1,h}(\Sigma_{f,0})\textrm{ and }\d_h \in [\chi_{2,h}(\Omega_{0})]^3$.

Here, the pair $(\u_h,p_h) \in [\chi_{2,h}(\Omega_0)]^3 \times \chi_{1,h}(\Omega_0)$ form the Taylor-Hood elements \cite{taylor_numerical_1973}, which are well known to be numerically stable for the incompressible Navier-Stokes problem while converging to the exact solution at the second order with respect to $h$ in their respective norms if the exact solution is sufficiently regular \cite{girault_finite_1986}.

As for $\lambda_h \in \chi_{1,h}(\Sigma_{f,0})$, although verification of the continuous inf-sup in the $(H^{1/2}(\Sigma_{f,0})'$ norm is a classical result by Verfürth \cite{Verfürth1986/87}, a proof of discrete stability in our specific case of curved 3D mesh elements was not found in the literature. On the other hand, Verfürth also proved in \cite{Verfürth1986/87} that a $[\chi_{1,h}(\Omega_0)]^3$ velocity field enriched by a single facet bubble function is stable when paired with a $\chi_{1,h}(\Sigma_{f,0})$ Lagrange multiplier $\lambda_h$. Since a velocity field in $[\chi_{2,h}(\Omega_0)]^3$ also possesses additional degrees of freedom along all facets of $\Sigma_{f,0}$, it is reasonable to conjecture that this coupling yields a stable finite element pair. This assumption, in addition to that of higher-order regularity, yields an order of convergence of 5/2 of the multiplier with respect to $h$ in the $(H^{1/2}(\Sigma_{f,0}))'$ norm (see Theorem 5.8.3 of \cite{brenner_mathematical_2008}).

On the other hand, $\d_h \in [\chi_{2,h}(\Omega_0)]^3$, unlike $p_h$ and $\lambda_h$, is not a Lagrange multiplier of $\u_h$, but is associated with a Laplace problem in the bulk of $\Omega_0$ and a Laplace-Beltrami operator over $\Sigma_{f,0}$, naturally providing $\d_h$ with the norm $\|.\|_M$. Rather than being required to satisfy an inf-sup condition, to guarantee existence and uniqueness, the Schur complement associated with $\d_h$ would have to be proven to be invertible, something that is outside of the scope of this article. Instead, hereafter, the linearized system of Equation (\ref{Eq:ParametrizedWeakForm}) is simply assumed to be well-posed. Further assuming sufficient solution regularity, $\d_h$ converges at the second order in the norm $\|.\|_M$.

Given these finite element spaces, the Jacobian of System (\ref{Eq:ParametrizedWeakForm}) may be discretized into the linear system below.

\begin{equation} \label{Eq:JacobianSystem}
    \underbrace{
    \begin{pmatrix}
        J_{uu} & J_{up} & J_{u\lambda} & J_{ud} & 0 \\
        J_{pu} & 0 & 0 & J_{pd} & 0 \\
        J_{\lambda u} & 0 & 0 & J_{\lambda d} & 0 \\
        0 & 0 & J_{d\lambda} & J_{dd} & J_{dp_a} \\
        0 & 0 & 0 & J_{p_a d} & J_{p_a p_a} \\
    \end{pmatrix}
    }_{=:J}
    \begin{pmatrix}
    \vec{\u}_h \\ \vec{p}_h \\ \vec{\lambda}_h \\ \vec{d}_h \\ p_a
    \end{pmatrix}
    =
    \begin{pmatrix}
        \vec{\boldsymbol{b_u}} \\ \vec{b_p} \\ \vec{b_\lambda} \\ \vec{\boldsymbol{b_d}} \\ b_{p_a}
    \end{pmatrix}
\end{equation}
where $J$ is the Jacobian matrix. In what follows, the subscript $h$ denotes a finite element approximation of characteristic mesh size $h$. Furthermore, the vector of nodal degrees of freedom of a finite element approximation is denoted by a vector arrow and subscript $h$, \textit{e.g.} $\vec{\u}_h$ is $\u_h$'s vector of coefficients. 

At each iteration of Newton's method, Linear System (\ref{Eq:JacobianSystem}) must be solved.

\section{Preconditioning} \label{Sec:Preconditioning}

In this section, after shortly summarizing Krylov-space methods and the role of preconditioning in their use, a preconditioner is derived for the solution of Equation (\ref{Eq:JacobianSystem}).

When solving a generic, sparse linear system ensuing from a PDE discretization such as 
\begin{equation}
    A\vec{z_h} = \vec{f},
\end{equation}
as defined in (\ref{Eq:GenericLinearSystem}), direct inversion is often the preferred method for problems of moderate size, \textit{i.e.} 2D or small 3D problems. At this size, direct methods are favored for their superior speed, exactness and because they can be implemented as black-box algorithms, agnostic to the PDE in question. Unfortunately, these benefits degrade at larger problem sizes where both the memory requirements and running time of direct approaches rapidly become prohibitive. 

Krylov-space methods, on the other hand, are a memory-efficient, iterative alternative at larger system sizes. Rather than relying on factorization algorithms to construct an exact solution, they employ a matrix-vector multiplication at each iteration to build up a subspace of $\mathbb{R}^N$, named the Krylov subspace. How the approximate solution $\vec{z_h}$ is formed from the subspace depends on which Krylov-space algorithm is employed. The iterative method is then considered to have converged once the residual has fallen below a user-prescribed tolerance.

Among Krylov-space methods are, for example, the Conjugate Gradient (CG) method, the Minimum Residual (MinRes) method, the Generalized Minimum Residual (GMRes) method and others. These methods yield convergent iterations under different circumstances. For instance, the CG algorithm is optimally convergent if $A$ is symmetric positive-definite, while the MinRes method is optimal for symmetric indefinite systems. On the other hand, GMRes, the method used here, is optimal for nonsymmetric systems. The convergence rate of these methods, provided that they are used under suitable conditions, is then governed by the conditioning of $A$. Convergence is fast when $A$ is well-conditioned, \textit{i.e.} when the minimal polynomial of the preconditioned operator is of low degree \cite{murphy_note_2000} and/or when the eigenvalues are clustered away from zero \cite{benzi_numerical_2005}. Conversely, convergence can be prohibitively slow when $A$ is ill-conditioned. The statement on the minimal polynomial is contained in an exact result for GMRes stating that, for a given iteration $k$ and current iterate $\vec{z}_{h,k}$, the residual $\vec{r}_{k} =A\vec{z}_{h,k} - \vec{f}_k$ satisfies
\begin{equation} \label{Eq:MinPolynomialGMRes}
    \|\vec{r}_k\|_2 = \underset{p_k\in\Pi_k,p_k(0) = 1}{\textrm{min}}\|p_k(A)\vec{r}_0\|_2
\end{equation}
where $\Pi_k$ is the space of polynomials of degree $k$ and $\|.\|_2$ denotes the Euclidean norm. This result can be found in classical references such as \cite{saad_gmres_1986,elman_finite_2014}.

It is well known that discretized PDEs are inherently ill-conditioned and, thus, do not lend themselves to Krylov-space solution in their native form. A widely-used approach, as is done here, is to solve an equivalent, "left-preconditioned" problem,
\begin{equation}
    P^{-1}A\vec{z_h} = P^{-1}\vec{f},
\end{equation}
where $P\in \mathbb{R}^{N\times N}$, the preconditioner, is a user-prescribed matrix intended to be (i) efficient to invert and (ii) such that $P^{-1}A$ is well-conditioned.

In the remainder of this section, a physics-based preconditioner is derived based on a block factorization of the Jacobian matrix in the high-capillarity limit. An approximate version of this preconditioner is then given together with an inversion algorithm.

\subsection{Physics-Based Preconditioner} \label{Sec:HighPreconditioner}

The block structure of the preconditioner must be adapted to the Krylov-space algorithm. The Jacobian $J$ in Equation (\ref{Eq:JacobianSystem}) is a nonsymmetric matrix, assumed to be invertible. An optimal algorithm, therefore, is GMRes. Furthermore, since there are no symmetry requirements on the system, a block triangular structure of the preconditioner is permissible. 

One may observe in Equation (\ref{Eq:ParametrizedWeakForm}) that in the high-capillarity limit, $J_{d\lambda}\overset{\Ca\rightarrow 0}{\rightarrow}0$. This suggests that in the limit $\Ca \ll 1$, $\d$ and $p_a$ are weakly coupled to $\u$, $p$ and $\lambda$. A preconditioner is now derived based on this hypothesis and, hence, on a block factorization of the Jacobian in this limit.
\begin{equation} \label{Eq:JacobianLimit}
    J \overset{\Ca\rightarrow 0}{\rightarrow} \begin{pmatrix}
        J_{uu} & J_{up} & J_{u\lambda} & J_{ud} & 0 \\
        J_{pu} & 0 & 0 & J_{pd} & 0 \\
        J_{\lambda u} & 0 & 0 & J_{\lambda d} & 0 \\
        0 & 0 & 0 & J_{dd} & J_{dp_a} \\
        0 & 0 & 0 & J_{p_a d} & J_{p_a p_a} \\
    \end{pmatrix}
\end{equation}
Thus, in this limit, $\d$ and $p_a$ are decoupled from $\u$, $p$ and $\lambda$. An appropriate approach to deriving the preconditioner is to treat the two sub-systems separately. Let us define the Navier-Stokes, free-surface and coupling blocks, respectively,
\begin{equation}
    N = \begin{pmatrix}
        J_{uu} & J_{pu}^T & J_{\lambda u}^T \\
        J_{pu} & 0 & 0 \\
        J_{\lambda u} & 0 & 0 \\
    \end{pmatrix} 
    \textrm{, }
    F = \begin{pmatrix}
        J_{dd} & J_{dp_a} \\
        J_{p_a d} & J_{p_a p_a} \\
    \end{pmatrix}
    \textrm{ and }
    K = \begin{pmatrix}
        J_{ud} & 0 \\
        J_{pd} & 0 \\
        J_{\lambda d} & 0
    \end{pmatrix}
\end{equation}
such that
\[
J \overset{\Ca\rightarrow0}{\rightarrow} \begin{pmatrix}
    N & K \\ 0 & F    
\end{pmatrix}.
\]

\subsubsection{Navier-Stokes Factorization}

$N$, here, is a double saddle-point block matrix. In order to derive the block preconditioner, $P_N$, for $N$, a Schur complement factorization is carried out. Let
\begin{align}
    A = J_{uu} 
    \textrm{, }
    B = \begin{pmatrix}
        J_{pu} \\
        J_{\lambda u}
    \end{pmatrix}
    \textrm{ and }
    C = 0.
\end{align}
Thus, 
\[
N = \begin{pmatrix}
    A & B^T \\
    B & C
\end{pmatrix},
\]
and a classical lower-triangular factorization of $N$ is
\begin{equation}
    N = \underbrace{\begin{pmatrix}
        A & 0 \\
        B & S
    \end{pmatrix}}_{=M}
    \begin{pmatrix}
        I & A^{-1}B^T \\
        0 & I
    \end{pmatrix},
\end{equation}
where $S$, the Schur complement, satisfies $S = C - BA^{-1}B^T$. Straightforward matrix algebra shows that \[(M^{-1}N - I)^2 = 0.\]
All eigenvalues of $M^{-1}N$ are, therefore, 1 and the minimal polynomial of $M^{-1}N$ is of the second order. This implies that GMRes applied to this operator would converge in exactly two iterations. 

Computing the inverse action of $M$ requires sequentially inverting matrix $A$ and matrix $S$. Unfortunately, $S$, now expanded,
\[
S = - \begin{pmatrix}
    J_{p u}J_{uu}^{-1}J_{p u}^T & J_{p u}J_{uu}^{-1}J_{\lambda u}^T \\ J_{\lambda u}J_{uu}^{-1}J_{p u}^T & J_{\lambda u}J_{uu}^{-1} J_{\lambda u}^T
\end{pmatrix},
\]
couples the variables $p$ and $\lambda$. This fact would prohibit a decoupled solution scheme when computing the inverse-action of $M$, compromising efficiency. In \cite{ali_beik_iterative_2018}, the authors showed that, in the case where $J_{uu}$ is symmetric positive-definite (when $\Rey\rightarrow 0$ here), for a matrix
\[P_S = - \begin{pmatrix}
    J_{p u}J_{uu}^{-1}J_{p u}^T & 0 \\ 0 & J_{\lambda u}J_{uu}^{-1} J_{\lambda u}^T
\end{pmatrix},
\]
all eigenvalues of the preconditioned matrix $P_S^{-1}S$ lie in the interval $(0,2)$. Furthermore, the eigenvalues are symmetric about $1$, satisfying $\lambda_i(P_S^{-1}S) = 1 \pm\sigma_i(X)$ where $\sigma_i(X)$ denotes the $i$-th singular value of matrix $X$ and 
\[X =  (J_{p u}J_{uu}^{-1}J_{p u}^T)^{-1/2} J_{pu}J_{uu}^{-1}J_{\lambda u}^T (J_{\lambda u}J_{uu}^{-1} J_{\lambda u}^T)^{-1/2}. \]
There is no lower bound on the minimum eigenvalue however, which may be arbitrarily close to zero without additional assumptions. The authors further explained in \cite{ali_beik_iterative_2018} that the minimum eigenvalue is close to 1 when the two Lagrange multipliers are weakly coupled. Given that the pressure constraint acts on bulk values of $\u$, while that of $\lambda$ acts on $\u$'s normal trace, it is reasonable to believe that these constraints act on substantially different modes of the velocity field, save for certain modes concentrated in a low-dimensional subspace.

At this stage, one may write
\begin{equation}
\begin{split}
    N &= \begin{pmatrix}
        A & 0 \\
        B & P_S P_S^{-1}S
    \end{pmatrix}
    \begin{pmatrix}
        I & A^{-1}B^T \\
        0 & I
    \end{pmatrix} \\
    &= \underbrace{\begin{pmatrix}
        A & 0 \\
        B & P_S
    \end{pmatrix}}_{=P_N}
    \begin{pmatrix}
        I & 0 \\
        0 & P_S^{-1}S
    \end{pmatrix}
    \begin{pmatrix}
        I & A^{-1}B^T \\
        0 & I
    \end{pmatrix},
\end{split}
\end{equation}
and, therefore,
\begin{equation}
    P_N^{-1}N = 
    \begin{pmatrix}
        I & A^{-1}B^T \\
        0 & P_S^{-1}S
    \end{pmatrix}.
\end{equation}
Hence, in the limit where $\Rey\rightarrow 0$, $P_N$ has identically $1$ eigenvalues in the velocity block and eigenvalues $\lambda_i(P_S^{-1}S) = 1 \pm\sigma_i(X)$ in the Lagrange multiplier blocks.

\subsubsection{Free-Surface Factorization}

Meanwhile, a standard Schur complement factorization of $F$, the $(\d,p_a)$ operator, shows that
\begin{equation}
    P_F = \begin{pmatrix}
        J_{dd} & 0 \\
        J_{p_a d} & J_{p_a p_a} - J_{p_a d}J_{dd}^{-1} J_{d p_a}
    \end{pmatrix}
\end{equation}
satisfies $(P_F^{-1}F - I)^2 = 0$. 

\subsubsection{Jacobian Factorization}

The factorization of the approximate Jacobian
\begin{equation}
    \begin{pmatrix}
        N & K \\
        0 & F \\
    \end{pmatrix}
    =
    \underbrace{\begin{pmatrix}
        P_N & K(P_F^{-1}F)^{-1} \\
        0 & P_F \\
    \end{pmatrix}}_{=:P}
    \begin{pmatrix}
        P_N^{-1}N & 0 \\
        0 & P_F^{-1}F \\
    \end{pmatrix}
\end{equation}
proves that in the limit where $\Ca \rightarrow 0$, $P^{-1}J$ has identically 1 eigenvalues in the $\u$, $\d$ and $p_a$ blocks, as well as eigenvalues in (0,2) in the $p$ and $\lambda$ blocks in the limit $\Rey \rightarrow 0$. Since only flows at low Reynolds numbers are considered here, the eigenvalues of $P_S^{-1}S$ are not expected to deviate significantly from (0,2). Preconditioner $P$, however, possesses dense blocks, as evidenced by its expanded representation below.
\begin{equation} \label{Eq:IdealHighPreconditioner}
    P = \begin{pmatrix}
        J_{uu} & 0 & 0 & J_{ud} & -J_{ud}J_{dd}^{-1}J_{dp_a} \\
        J_{pu} & -S_p & 0 & J_{pd} & -J_{pd}J_{dd}^{-1}J_{dp_a} \\
        J_{\lambda u} & 0 & -S_\lambda  & J_{\lambda d} & - J_{\lambda d}J_{dd}^{-1}J_{dp_a} \\
        0 & 0 & 0 & J_{dd} & J_{dp_a} \\
        0 & 0 & 0 & 0 & J_{p_a p_a} - J_{p_a d}J_{dd}^{-1}J_{d p_a} \\
    \end{pmatrix},
\end{equation}
where 
\begin{equation} \label{Eq:SchurComplements}
    S_p := J_{pu}J_{uu}^{-1}J_{pu}^T\textrm{ and }S_\lambda := J_{\lambda u}J_{uu}^{-1}J_{\lambda u}^T.
\end{equation}
Indeed, matrix blocks $(2,2)$ and $(3,3)$ are dense due to the presence of $J_{uu}^{-1}$. Thus, $P$, an ideal preconditioner, may not realistically be implemented in scalable computations. Instead, approximations of these dense blocks by sparse counterparts must be made.

\subsection{Approximation} \label{Sec:PreconditionerApproximation}

In preconditioner $P$ as defined in Equation (\ref{Eq:IdealHighPreconditioner}), the dense matrices $S_p$, $S_\lambda$, must be approximated by sparse counterparts that are, in the best case, spectrally equivalent operators. 

\begin{definition} \label{Def:SpectralEquivalence}
\textbf{Spectral Equivalence} \cite{elman_finite_2014} \\
    A symmetric positive-definite (SPD) matrix $Y\in \mathbb{R}^{N\times N}$ is spectrally equivalent to an SPD $X\in \mathbb{R}^{N\times N}$ if there exist two positive constants $c_1$ and $c_2$, both independent of the mesh size such that
\begin{equation}
    c_1 \leq \frac{\vec{v}^T X\vec{v}}{\vec{v}^T Y\vec{v}} \leq c_2 \,\forall \vec{v}\in \mathbb{R}^N.
\end{equation}
\end{definition}
When matrix blocks in the preconditioner are replaced by spectrally equivalent approximations, the number of Krylov-space iterations is expected to be independent of the mesh size. 

Operator $S_p$, as defined in Equation (\ref{Eq:SchurComplements}), is the pressure Schur complement that arises in the Navier-Stokes equations. It is well known that in the limit $\Rey \rightarrow 0$, $S_p$ is spectrally equivalent to the pressure mass matrix, $M_p$ \cite{elman_finite_2014}. This statement is a direct consequence of the hypotheses of Theorem \ref{Thm:WellPosedness} applied to $(\u_h,p_h)$, together with the fact that the pressure mass matrix induces its natural norm, \textit{i.e.} \[\|p_h\|_Q = \sqrt{\vec{p}_h^TM_p \vec{p}_h},\] and that $J_{uu}$ induces the velocity norm, \textit{i.e.} \[ \|\u_h\|_{V_D} = \sqrt{\vec{\u}_h^T J_{uu}\vec{\u}_h}.\]

\begin{lemma} \label{Lemma:SpectralEquivalence}
    Let $V_h \times Q_h$ denote the Taylor-Hood elements, and let $J_{pu}$ and $J_{uu}$ be as defined in Equation (\ref{Eq:JacobianSystem}) in the limit where $\Rey \rightarrow 0$. Then, the Schur complement, $J_{pu}J_{uu}^{–1}J_{pu}^T$ is spectrally equivalent to the pressure mass matrix, $M_p$, such that
    \[
    \beta_{J_{pu}}^2 \leq \frac{\vec{p_h}^T J_{pu}J_{uu}^{–1}J_{pu}^T\vec{p_h}}{\vec{p_h}^T M_p \vec{p_h}} \leq \gamma_{J_{pu}}^2 \,\forall p_h\in Q_h,
    \]
    where $\gamma_{J_{pu}}$ and $\beta_{J_{pu}}$ are the continuity and inf-sup constants, respectively, associated with the operator $J_{pu}$.
    \begin{proof}
    Recalling that the Taylor-Hood elements satisfy Theorem \ref{Thm:WellPosedness}, for any $\vec{p_h} \in Q_h$,
        \begin{equation}
    \begin{split}
        \frac{\vec{p}_h^T J_{pu}J_{uu}^{-1}J_{pu}^T \vec{p}_h}{\|p_h\|_{Q}^2} &\leq \underset{\vec{p}_h\in \mathbb{R}^{\textrm{dim}(Q_h)}}{\textrm{sup}} \frac{\vec{p}_h^T J_{pu}J_{uu}^{-1}J_{pu}^T \vec{p}_h}{\|p_h\|_{Q}^2} \\
        &\overset{J_{uu}\textrm{ is SPD}}{=} \underset{\vec{p}_h\in \mathbb{R}^{\textrm{dim}(Q_h)}}{\textrm{sup}} \frac{(\vec{p}_h^T J_{pu}J_{pu}^T \vec{p}_h)^2}{(\vec{p}_h^T J_{pu}J_{uu}J_{pu}^T \vec{p}_h)\|p_h\|_{Q}^2} \\
        &=\underset{ \vec{p}_h\in \mathbb{R}^{\textrm{dim}(Q_h)} }{\textrm{sup}}\underset{ \vec{w}_h\in \mathbb{R}^{\textrm{dim}(V_h)} }{\textrm{sup}}\frac{(\vec{w}_h^T J_{pu}^T \vec{p}_h)^2}{ (\vec{w}_h^T J_{uu} \vec{w}_h)  \|p_h\|_{Q}^2} \\
        &=\underset{ \vec{p}_h\in \mathbb{R}^{\textrm{dim}(Q_h)} }{\textrm{sup}}\underset{ \vec{w}_h\in \mathbb{R}^{\textrm{dim}(V_h)} }{\textrm{sup}}\frac{(-\int_{\Omega_0}(\bnabla\cdot w_h)p_h dV )^2}{ \|w_h\|_{V_D}^2  \|p_h\|_{Q}^2} \overset{\textrm{Continuity}}{\leq} \gamma_{J_{pu}}^2. 
    \end{split}
\end{equation}
Operator $J_{uu}$ is SPD, as it is the matrix induced by a symmetric, coercive bilinear form. For any $p_h\in Q_h$,
\begin{equation}
    \begin{split}
        \frac{\vec{p}_h^T J_{pu}J_{uu}^{-1}J_{pu}^T \vec{p}_h}{\|p_h\|_{Q}^2} &\geq \underset{\vec{p}_h\in \mathbb{R}^{\textrm{dim}(Q_h)}}{\textrm{inf}} \frac{\vec{p}_h^T J_{pu}J_{uu}^{-1}J_{pu}^T \vec{p}_h}{\|p_h\|_{Q}^2} \\
        &\overset{J_{uu}\textrm{ is SPD}}{=} \underset{\vec{p}_h\in \mathbb{R}^{\textrm{dim}(Q_h)}}{\textrm{inf}} \frac{(\vec{p}_h^T J_{pu}J_{pu}^T \vec{p}_h)^2}{(\vec{p}_h^T J_{pu}J_{uu}J_{pu}^T \vec{p}_h) \|p_h\|_{Q}^2} \\
        &= \underset{ \vec{p}_h\in \mathbb{R}^{\textrm{dim}(Q_h)} }{\textrm{inf}} \underset{\vec{w}_h\in \mathbb{R}^{\textrm{dim}(V_h)}}{\textrm{sup}} \frac{(\vec{w}_h^T J_{pu}^T \vec{p}_h)^2}{(\vec{w}_h^TJ_{uu}\vec{w}_h) \|p_h\|_{Q}^2} \\
        &= \underset{\vec{p}_h\in \mathbb{R}^{\textrm{dim}(Q_h)}}{\textrm{inf}} \underset{\vec{w}_h\in \mathbb{R}^{\textrm{dim}(V_h)} }{\textrm{sup}} \frac{(-\int_{\Omega_0}(\bnabla\cdot w_h)p_h dV )^2}{\|w_h\|_{V_D}^2 \|p_h\|_{Q}^2} \overset{\textrm{inf-sup}}{\geq} \beta_{J_{pu}}^2.
    \end{split}
\end{equation}  
    \end{proof}
\end{lemma}

Hence, $S_p$ is spectrally equivalent to the operator inducing the pressure norm, which in this case, is the pressure mass matrix. Considering that $\textrm{diag}(M_p)$ is, itself, spectrally equivalent to $M_p$ \cite{elman_finite_2014} and faster to invert, the approximation
\[\hat{S}_p = \textrm{diag}(M_p)\]
is implemented to approximate $S_p$. 

One may observe that, akin to $S_p$, $S_\lambda$ is also a Schur complement with the same structure, whose constituents are assumed here to satisfy the hypotheses of Theorem \ref{Thm:WellPosedness}. The same proof as in Lemma \ref{Lemma:SpectralEquivalence} may therefore be carried out to derive a spectrally equivalent approximation of $S_\lambda$, giving
\begin{equation}
    \beta_{J_{\lambda u}}^2 \leq \frac{\vec{\lambda}_h^T S_\lambda \vec{\lambda}_h }{\|\lambda_h \|^2_{L_0}} \leq \gamma_{J_{\lambda u}}^2\,\forall \lambda_h \in L_{0,h}.
\end{equation}
Therefore, in the limit $\Rey\rightarrow 0$, a spectrally equivalent operator of $S_\lambda$, $\hat{S}_{\lambda}$, is one satisfying
\begin{equation}
\begin{split}
    \vec{\lambda}_h^T \hat{S}_\lambda \vec{\lambda}_h &= \|\lambda_h\|^2_{L_0} \\
    &= \langle (-\Delta_{S,0})^{-1/4} \lambda_h, (-\Delta_{S,0})^{-1/4} \lambda_h \rangle_{L^2(\Sigma_{f,0})} \\
    &= \langle \lambda_h, (-\Delta_{S,0})^{-1/2} \lambda_h \rangle_{L^2(\Sigma_{f,0})},
\end{split}
\end{equation}
where in the final step, the self-adjointness of the Laplace-Beltrami operator is used. Introducing the finite element mass and stiffness matrices for $\lambda_h$ as $M_\lambda$ and $K_\lambda$, respectively, the discrete norm may be expressed as
\begin{equation}
\begin{split}
    \|\lambda_h\|^2_{L_0} &= \vec{\lambda}_h^T  M_\lambda^{1/2} (M_\lambda^{-1/2}K_\lambda M_\lambda^{-1/2})^{-1/2} M_\lambda^{1/2} \vec{\lambda}_h.
\end{split}
\end{equation}
By definition of the $L_0$ norm in Section \ref{Sec:NavierStokes}, $K_\lambda$ has homogeneous Dirichlet boundary conditions along $\Gamma_D$ and homogeneous Neumann boundary conditions along $\partial \Sigma_f \setminus\Gamma_D$. Hence,
\begin{equation}
    \hat{S_\lambda} = M_\lambda^{1/2} ( M_\lambda^{-1/2}K_\lambda M_\lambda^{-1/2} )^{-1/2} M_\lambda^{1/2}
\end{equation}
is a spectrally equivalent approximation of $S_\lambda$ and the inverse action of $\hat{S}_\lambda$ that must be effected at each iteration of GMRes is
\begin{equation}
    \hat{S}_\lambda^{-1} = M_\lambda^{-1/2} ( M_\lambda^{-1/2}K_\lambda M_\lambda^{-1/2})^{1/2} M_\lambda^{-1/2}.
\end{equation}
In practice, this matrix is not explicitly assembled, and only its action on a vector is computed. How this action is computed is presented in Section \ref{Sec:Slambda}.

In the $p_a$ column of $P$, although $J_{dd}^{-1}$ is a dense operator, it acts on a column matrix $J_{d p_a}$. Thus, by approximately solving the system 
\[J_{dd}X = J_{d p_a}\] 
via one multi-grid V-cycle during the preconditioner setup phase, the ensuing column matrix, $X$ may be computed just once.

These choices result in the approximate preconditioner found below.

\begin{equation} \label{Eq:ApproximateHighPreconditioner}
    \hat{P} = \begin{pmatrix}
        J_{uu} & 0 & 0 & J_{ud} & -J_{ud}J_{dd}^{-1}J_{dp_a} \\
        J_{pu} & -\hat{S_p} & 0 & J_{pd} & - J_{pd}J_{dd}^{-1}J_{dp_a} \\
        J_{\lambda u} & 0 & -\hat{S_\lambda} & J_{\lambda d} & - J_{\lambda d}J_{dd}^{-1}J_{dp_a} \\
        0 & 0 & 0 & J_{dd} & J_{dp_a} \\
        0 & 0 & 0 & 0 & J_{p_a p_a} - J_{p_a d}J_{dd}^{-1}J_{d p_a} \\
    \end{pmatrix}
\end{equation}

\subsection{Inversion}

To ensure a swift application of preconditioner inverse $\hat{P}^{-1}$, the following sequential inversion algorithm is effected.

\begin{equation} \label{Eq:PreconditionInversion}
    \begin{split}
        &\underline{\textrm{\textbf{Inversion of $\hat{P}$}}} \\
        &\begin{cases}
            \textrm{Scalar inversion of } (J_{p_a p_a} - J_{p_a d}J_{dd}^{-1}J_{dp_a}) p_a = b_{p_a} \\
            \\
            \textrm{Sequential solve of } J_{dd}\vec{\d}_h = \vec{\boldsymbol{b_d}} - J_{dp_a}p_a \\
            \\
            \textrm{Multi-grid solve of }J_{uu}\vec{\u}_h = \vec{\boldsymbol{b_u}} - J_{ud}\vec{\d}_h + J_{ud}J_{dd}^{-1}J_{d p_a}p_a \\
            \\
            \textrm{Diagonal solve of }-\textrm{diag}(M_p)\vec{p}_h = \vec{b_p} - J_{pu}\vec{\u}_h \\
            \\
            \textrm{Rational approximation of } - \hat{S}_\lambda \vec{\boldsymbol{\lambda}_h} = \vec{b_\lambda} - J_{\lambda d}\vec{\d}_h + J_{\lambda d}J_{dd}^{-1}J_{d p_a}p_a - J_{\lambda u}\vec{\u}_h \\
        \end{cases}
    \end{split}
\end{equation}
Multi-grid inversions of $J_{uu}$ are performed with PETSc's implementation of the Hypre BoomerAMG algebraic multi-grid algorithm. Specifically, one V-cycle is implemented with a Chebyshev smoother together with HMIS coarsening and extended interpolation as the approximate solve. These choices are made in order to enhance the parallel scalability of the preconditioners. The sequential solve of $J_{dd}$ and inversion action of $\hat{S_\lambda}$ require special treatment and are explained next.


\subsubsection{Approximate Displacement Solve} \label{Sec:DisplacementSolve}

Recalling that $\d_h$ satisfies a Laplace-Beltrami equation over $\Sigma_{f,0}$ for tangent components and a Young-Laplace equation for the normal component, $\d_h$'s degrees of freedom over $\Sigma_{f,0}$ are therefore decoupled from interior values in operator $J_{dd}$ of (\ref{Eq:ApproximateHighPreconditioner}), leading to the block structure
\[
    J_{dd}=\begin{pmatrix}
        J_{\Omega_0 \Omega_0} & J_{\Omega_0 \Sigma_{f,0}} \\
        0 & J_{\Sigma_{f,0} \Sigma_{f,0}}
    \end{pmatrix}.
\]
Monolithic multi-grid algorithms may fail to produce mesh-independent iterations due to disparity of bulk and free-surface operators. Instead, within a preconditioner application, $J_{dd}$ may be inverted via a sequential multi-grid algorithm.
\begin{equation}
    \begin{split}
        &\underline{\textrm{\textbf{Inversion of $J_{dd}$}}} \\
        &\begin{cases}
             \textrm{Multi-grid solve of }J_{\Sigma_{f,0} \Sigma_{f,0}}\vec{\d}_{h,\Sigma_{f,0}} =  \vec{\boldsymbol{b}}_{\d,\Sigma_{f,0}} \\
            \\
            \textrm{Multi-grid solve of }  J_{\Omega_{0} \Omega_{0}}\vec{\d}_{h,\Omega_{0}} =  \vec{\boldsymbol{b}}_{\d,\Omega_{0}} - J_{\Omega_0 \Sigma_{f,0}}\vec{\d}_{h,\Omega_{0}} \\
        \end{cases}
    \end{split}
\end{equation}
Multi-grid inversions are carried out with the same parameters as $J_{uu}$.

\subsubsection{Rational Approximation of a Matrix Square Root} \label{Sec:Slambda}

This section is an adaptation of the method presented in Section 4 of \cite{alexandrou_multigrid_2019}.

The inverse action of matrix $\hat{S_\lambda}$ entails computing the action of
\begin{equation}
    \vec{\lambda}_h = M_\lambda^{-1/2} (M_\lambda^{-1/2}K_\lambda M_\lambda^{-1/2})^{1/2} M_\lambda^{-1/2} \vec{\v}_\lambda
\end{equation}
for an input vector $\vec{\v}_\lambda$ at every iteration of GMRes. To lighten the notation in subsequent steps, let 
\begin{equation} \label{Eq:A}
    A = M_\lambda^{-1/2}K_\lambda M_\lambda^{-1/2},
\end{equation} 
giving rise to the inverse action
\begin{equation}
    \vec{\lambda}_h = M_\lambda^{-1/2} A^{1/2} M_\lambda^{-1/2} \vec{\v}_\lambda.
\end{equation}
Let us note that, by definition, $A$ is SPD. 

In what follows, an efficient approximation of the square-root action by operator $A$ is presented. Square-root functions may be accurately and efficiently approximated by a rational function $R_{m,n}$ of user-defined degrees $m$ and $n$ such that 
\begin{equation}
    R_{m,n}(z) = C \frac{\prod_{i=1}^m (z+n_i) }{\prod_{j=1}^n (z+d_j)},
\end{equation}
where $C$, $n_i$ and $d_i$ are constants for $i=1,...,n$. In the case where the numerator and denominator degrees are equal, \textit{i.e.} $m=n$, Zolotarev derived an analytical solution in \cite{zolotarev_egor_ivanovich_application_1877} (see \cite{higham_functions_2008} for a more recent presentation of this result),
\begin{equation}
\begin{split}
    R_{n,\epsilon}(z) :&= \underset{R_{n,n}}{\textrm{argmin}}\underset{\Tilde{z}\in(\epsilon,1)}{\textrm{max}} |1 - \sqrt{\Tilde{z}} R_{n,n}(\Tilde{z}) | \\
    &= a_{n,\epsilon} \prod_{j=1}^n \frac{z + c_{n,\epsilon,(2j-1)}}{z+c_{n,\epsilon,2j}},
\end{split}
\end{equation}
that minimizes the $L^\infty$ error of the inverse square root function over an interval $(\epsilon,1)$ for $\epsilon>0$. Here, the coefficients, as presented in \cite{alexandrou_multigrid_2019}, are given by
\begin{equation}
    c_{n,\epsilon,k} = \textrm{cs}^2(k\,v_{n,\epsilon},\sqrt{1-\epsilon})\textrm{ for $k=1,...,2n$}
\end{equation}
\[\textrm{with }v_{n,\epsilon} = \frac{K(\sqrt{1-\epsilon})}{2n+1},\]
\begin{equation}
    a_{n,\epsilon} = \frac{2}{1+\sqrt{1-d_{n,\epsilon}^2}} \prod_{j=1}^n \frac{s_{n,\epsilon,(2j-1)}}{s_{n,\epsilon,2j}}
\end{equation}
\[\textrm{with }s_{n,\epsilon,k} = \textrm{sn}^2(k\,v_{n,\epsilon},\sqrt{1-\epsilon}) \textrm{for $k=1,...,2n$ and }d_{n,\epsilon} = (1-\epsilon)^{\frac{2n+1}{2}}\prod_{j=1}^n s_{n,\epsilon,(2j-1)}^2.\]

Here, $\textrm{sn}(u,k)$ and $\textrm{cs}(u,k)  = \textrm{cn}(u,k)/\textrm{sn}(u,k)$ are the Jacobi elliptic functions and $K(k)$ is the complete elliptic integral of the first kind. Given their dependence on parameter $\epsilon$, these quantities are spectrum-dependent when applied to a matrix argument. They must therefore be computed once during preconditioner assembly. In simulations presented here, they are calculated with built-in routines from the SciPy library \cite{virtanen_scipy_2020}. Furthermore, an analytical expression of the maximum point-wise error exists,
\begin{equation}
    \underset{\epsilon < z < 1}{\textrm{max}}|1-\sqrt{z}R_{n,\epsilon}(z)| = \frac{d_{n,\epsilon}^2}{1+\sqrt{1-d_{n,\epsilon}^2}}
\end{equation}
and may also be computed at assembly. The maximum $L^\infty(\epsilon,1)$ error therefore decays exponentially with respect to the degree of the approximation. Via a partial fraction expansion, $R_{n,\epsilon}$ may then equivalently be expressed as
\begin{equation}
    R_{n,\epsilon}(z) = a_{n,\epsilon}( 1 + \sum_{i=1}^n\frac{q_i}{z-p_i} )
\end{equation}
\[\textrm{where } q_i = (c_{n,\epsilon,(2i-1)}-c_{n,\epsilon,2i})\prod_{j=1,j\neq i}^n \frac{c_{n,\epsilon,(2j-1)} - c_{n,\epsilon,2i}}{c_{n,\epsilon,2j} - c_{n,\epsilon,2i}} \textrm{ and }p_i  = - c_{n,\epsilon,2i}.\]
Let $\lambda_\textrm{min}$ and $\lambda_\textrm{max}$ denote the minimum and maximum eigenvalues, respectively of the SPD matrix $A$. Defining the matrix 
\[\Tilde{A} = \frac{1}{\lambda_\textrm{max}}A,\]
with an eigenvalue spectrum $(\frac{\lambda_\textrm{min}}{\lambda_\textrm{max}},1)$, Zolotarev's solution for the inverse square root is therefore optimal when applied to the argument $\Tilde{A}$ and minimizes the rational approximation's $L^\infty$ error over its spectrum. The square root action of $A$ may then be approximated via the relation,
\begin{equation}
    \begin{split}
        A^{1/2} &= A A^{-1/2} \\
        &= \frac{1}{\sqrt{\lambda_\textrm{max}}} A (\frac{1}{\lambda_\textrm{max}}A)^{-1/2} \\
        &\approx \frac{1}{\sqrt{\lambda_\textrm{max}}} A R_{n,\epsilon}(\frac{1}{\lambda_\textrm{max}}A) \\
        &= \frac{a_{n,\epsilon}}{\sqrt{\lambda_\textrm{max}}}A\{ I + \sum_{i=1}^n q_i (\frac{1}{\lambda_\textrm{max}}A - p_iI)^{-1}  \} \\
        &= \frac{a_{n,\epsilon}}{\sqrt{\lambda_\textrm{max}}}M_\lambda^{-1/2}K_\lambda M_\lambda^{-1/2}\{ I + \sum_{i=1}^n q_i(\frac{1}{\lambda_\textrm{max}}M_\lambda^{-1/2}K_\lambda M_\lambda^{-1/2} - p_iI)^{-1}  \} \\
        &= \frac{a_{n,\epsilon}}{\sqrt{\lambda_\textrm{max}}}M_\lambda^{-1/2}K_\lambda M_\lambda^{-1/2}\{ I + \sum_{i=1}^n q_i M_\lambda^{1/2}(\frac{1}{\lambda_\textrm{max}} K_\lambda - p_i M_\lambda )^{-1}M_\lambda^{1/2}   \} \\
    \end{split}
\end{equation}
where $\epsilon = \lambda_\textrm{min}/\lambda_\textrm{max}$. The full inverse action of $\hat{S_\lambda}$ then becomes
\begin{equation}
    \begin{split}
        \hat{S_\lambda}^{-1} &= M_\lambda^{-1/2} A^{1/2} M_\lambda^{-1/2} \\
        &= M_\lambda^{-1/2} \frac{a_{n,\epsilon}}{\sqrt{\lambda_\textrm{max}}}M_\lambda^{-1/2} K_\lambda M_\lambda^{-1/2}\{ I + \sum_{i=1}^n q_i M_\lambda^{1/2}(\frac{1}{\lambda_\textrm{max}}K_\lambda - p_i M_\lambda )^{-1}M_\lambda^{1/2}   \} M_\lambda^{-1/2} \\
        &= \frac{a_{n,\epsilon}}{\sqrt{\lambda_\textrm{max}}}M_\lambda^{-1} K_\lambda \{ M_\lambda^{-1} + \sum_{i=1}^n q_i \underbrace{(\frac{1}{\lambda_\textrm{max}} K_\lambda + |p_i| M_\lambda}_{=:T_i} )^{-1} \}.
    \end{split}
\end{equation}
The expression is devoid of any fractional powers and may be computed via a sum of matrix-vector solves and a matrix-vector product. Each matrix $T_i$ is SPD and approximately inverted by the CG algorithm preconditioned by one V-cycle of PETSc's Hypre BoomerAMG algorithm. The mass matrix $M_\lambda$ is also SPD and inverted by the CG algorithm with Jacobi preconditioning and a relative tolerance of $10^{-10}$. The low tolerance ensures that the preconditioner is approximately identical at each iteration of GMRes. A degree $n$ of 10 is used in numerical experiments here for the rational approximation. Since $\lambda_h$ is a scalar linear finite element defined only over the free surface, approximately solving the SPD linear systems associated with $T_i$ and $M_\lambda$ is inexpensive compared to the velocity and displacement multi-grid cycles. 

Let us remark that $\hat{S}_\lambda$ is furthermore SPD and, therefore, may be used within Krylov-space methods for symmetric problems.

\section{Numerical Simulations} \label{Sec:NumericalSimulations}

In this section, the preconditioner derived in Section \ref{Sec:Preconditioning} is tested on two numerical examples of incompressible free-surface flows. The first one, the die-swell problem, is a numerical benchmark where a liquid exits a pipe into a quiescent atmosphere. The jet profile either expands or contracts depending on the flow rate and on liquid parameters.

In the second numerical example, a coating flow is modeled where a liquid exits a vertical pipe and is deposited onto a substrate moving in the horizontal direction. The liquid wets the substrate, forming a triple line between the solid, the liquid and surrounding gas. An equilibrium contact angle, dependent only upon the interfacial energies of the liquid, solid, and gas phases, is imposed as a parameter of the flow simulations.

Both numerical examples involve a stick-slip transition where the liquid undergoes a transition from a no-slip condition to a shear-free boundary condition. In the die-swell problem, this transition is known to engender integrable singularities both in velocity and pressure \cite{sprittles_finite_2012}. In the coating-flow example, however, the stick-slip phenomenon has more deleterious consequences on the solution, leading to a non-integrable singularity. This occurs as a result of the contact line paradoxically spreading while the liquid is simultaneously required to satisfy the no-slip condition, also known as the Huh-Scriven paradox \cite{huh_hydrodynamic_1971}. This non-integrable singularity is avoided, here, by means of a Navier-slip boundary condition in the vicinity of the contact line, one that relaxes the no-slip condition, causing a gradual transition from no slip to zero shear. In this case, Dussan and Davis proved that this leads to an integrable singularity \cite{dussan_v_motion_1974}.

So as to mitigate the discretization errors caused by the reduced regularity of the solution in the contact-line neighborhood, we choose to implement meshes that are finer in this region, but without having conducted a rigorous study of localized mesh convergence. This choice is made since the main focus of this article is preconditioning. 

\subsection{Implementation}

All spatial discretizations of the domain are performed with the library Gmsh \cite{geuzaine_gmsh_2009}. The Mmg library is implemented to optimize meshes that have become severely distorted throughout numerical continuation \cite{dapogny_three-dimensional_2014}. In parallel computations, the meshes are distributed across the different processes by the partitioner PT-Scotch. The weak forms are discretized by the FEM with FEniCSx into PETSc data structures \cite{logg_dolfin_2010,dalcin_parallel_2011}. Linear systems are solved with the preconditioners reported in Section \ref{Sec:Preconditioning}, which are all implemented in PETSc. The message-passing interface (MPI) is used for parallelized computation. All simulations are run on a high-performance computing platform.

\subsection{Die Swell} \label{Sec:DieSwell}

\begin{figure}
    \centering
    \includegraphics[width=0.6\linewidth]{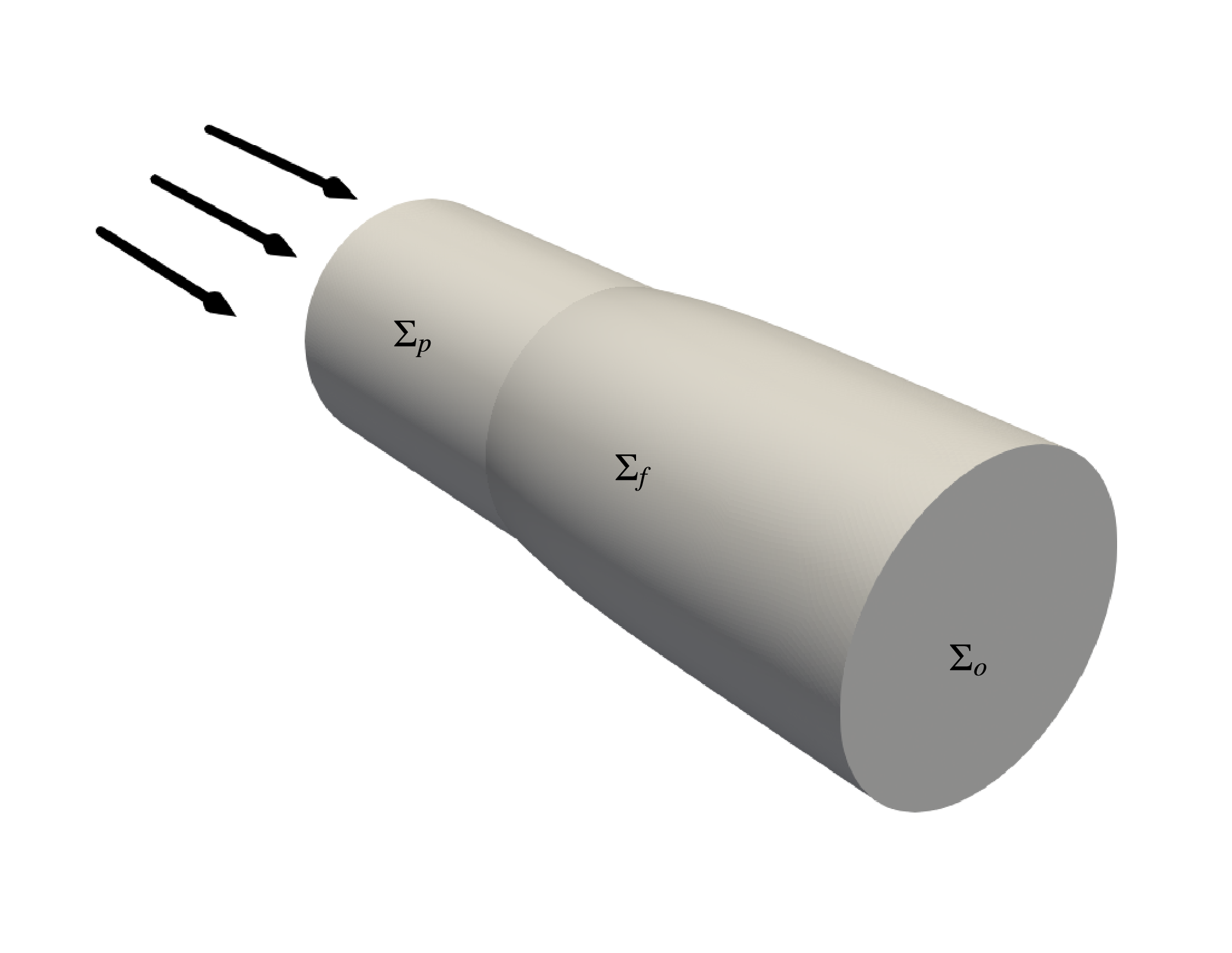}
    \caption{Die-swell geometry: after entering the pipe surface, $\Sigma_p$, the liquid exits into air and forms a free interface between air and liquid, $\Sigma_f$, before exiting through the numerical outlet, $\Sigma_o$. This axisymmetric problem is modeled in Cartesian coordinates without axisymmetry imposed.}
    \label{Fig:DieSwell}
\end{figure}

In the die-swell problem, after entering a pipe, $\Sigma_p$, with a Poiseuille velocity profile, a liquid jet is formed when an incompressible liquid exits the channel into a quiescent atmosphere. The free interface formed between the liquid and the gas, $\Sigma_f$, is governed by the Young-Laplace equation, a shear-free condition and an impermeability condition. The liquid then exits the numerical outlet, $\Sigma_o$, subject to a homogeneous Neumann boundary condition on the normal component and zero tangent velocity. For further details on this benchmark, see \cite{benkley_generalization_2026}.

In order to test the performance of the numerical method and preconditioner described respectively in Sections \ref{Sec:NumericalMethod} and \ref{Sec:Preconditioning}, axisymmetry is not exploited to simplify this problem which, instead, is solved in 3D, Cartesian coordinates.

Let us introduce the velocity spaces
\[\begin{split}
    V_{D,h} = \{ \v\in [\chi_{2,h}(\Omega_0)]^3 \mid\, &\v|_{\Sigma_{p,0}} = 2(1-x^2-y^2) \begin{pmatrix}0 \\ 0 \\ 1 \end{pmatrix}, \, v_x|_{\Sigma_{o,0}}= v_y|_{\Sigma_{o,0}} = 0\}
\end{split}\]
and 
\[\begin{split}
    V_{0,h} = \{ \v\in [\chi_{2,h}(\Omega_0)]^3 \mid\, &\v|_{\Sigma_{p,0}} = 0,\,v_x|_{\Sigma_{o,0}}= v_y|_{\Sigma_{o,0}} = 0\}.
\end{split}\]
Moreover, the displacement spaces are defined,
\[\begin{split}
    M_h = \{ \d\in [\chi_{2,h}(\Omega_0)]^3 \mid \, d_z|_{\Omega_0} = 0,\, &\d|_{\Sigma_{p,0}} = 0 \},
\end{split}\]
and
\[\begin{split}
    M_{\textrm{bulk},h} = \{ \d\in M_h \mid\,\d|_{\Sigma_{f,0}} = 0\}.
\end{split}\]
In this particular problem, deformation of the domain along the $z$ axis is not required to generate the space of solutions and is therefore made to be zero everywhere in $\Omega_0$, \textit{i.e.} $d_z|_{\Omega_0} = 0$. A homogeneous Neumann boundary condition is imposed on $\d$ over $\partial \Sigma_{f,0}\cap \partial \Sigma_{o,0}$ in the closure equation, \textit{i.e.}
\[ \t\cdot(\s_f\cdot\bnabla_S)\d = 0,\textrm{ along $\partial\Sigma_{f,0}\cap \partial \Sigma_{o,0}$}. \]
The closure equation, thus, is
\[
\int_{\Sigma_{f,0}} \bnabla_S\d_h:\bnabla_S[(I - \frac{\F_0^{-T}\n_0\n_0^T\F_0^{-1}}{|\F_0^{-T}\n_0|^2})\boldsymbol{\xi}_h] dS = 0,\,\forall \xi_h \in M_{h}.
\]

The Weak Form (\ref{Eq:DieSwellWeakForm}), found below, is parametrized by $\Omega_0$ and then solved by the numerical method described in Section \ref{Sec:NumericalMethod}. 

\begin{equation} \label{Eq:DieSwellWeakForm}
    \begin{split}
    &\text{Find $(\u_h,p_h,\lambda_h,\d_h,p_a)\in V_{D,h}\times Q_h \times L_{0,h} \times M_h \times \mathbb{R} $ such that} \\
    &\begin{cases}
        \int_\Omega \{\,\Rey\,\v_h\cdot (\u_h\cdot \bnabla )\u_h  + [\bnabla\u_h+(\bnabla\u_h)^T]:\bnabla\v_h  -p_h(\bnabla\cdot\v_h)\}dV \\
        \text{\textcolor{white}{aaa}}- \int_{\Sigma_f} \lambda (\v_h\cdot\n) dS = 0,\,\hfill\forall \v_h\in V_{0,h} \\
        \int_\Omega -q_h(\bnabla \cdot \u_h)dV = 0,\,\hfill\forall q_h\in Q_h \\
        \int_{\Sigma_f} -\chi_h(\u_h\cdot \n)dS = 0,\,\hfill\forall \chi_h\in L_{0,h} \\
        \int_{\Sigma_{f,0}} \bnabla_S\d_h:\bnabla_S[\boldsymbol{\xi}_h-(\boldsymbol{\xi}_h\cdot \frac{\F_0^{-T}\n}{|\F_0^{-T}\n|}) \frac{\F_0^{-T}\n}{|\F_0^{-T}\n|}]dS \\
        \text{\textcolor{white}{aaa}}+\int_{\Sigma_f} [ (\boldsymbol{\xi}_h\cdot\n)( \Ca \,\lambda_h + p_a) + \bnabla_S\cdot\boldsymbol{\xi}_h  ] dS = 0,\,\hfill\forall \boldsymbol{\xi}_h\in M_{h} \\
        \int_{\Omega_0} \bnabla\d_h:\bnabla\boldsymbol{\xi}_h\,dV = 0,\,\hfill\forall \boldsymbol{\xi}_h\in M_{\textrm{bulk},h} \\
        \int_{\partial \Sigma_o} p_a \,dl + 2\pi = 0
    \end{cases}
    \end{split}
\end{equation}
All Newton solves use the Stokes flow solution as the initial guess, where the initial free surface is an undeformed cylinder. The Young-Laplace condition is relaxed and the displacement field is omitted from starting point computations. None of the simulations reported here required numerical continuation when using this initial guess. Nonlinear iterations are considered to have converged when the maximum absolute value of the prescribed Newton correction falls below an absolute tolerance of $10^{-5}$ for $\u_h$ and $\d_h$. Krylov-space iterations are considered to have converged when the Euclidean norm of the residual falls below either a relative tolerance of $10^{-6}$ or an absolute tolerance of $10^{-13}$. Computations are carried out at a fixed Ohnesorge number, $Oh = \sqrt{\Ca/\Rey} = \sqrt{3}$ for comparison with the values reported in \cite{benkley_generalization_2026}. 

\begin{figure}
    \centering
    \includegraphics[width=0.6\linewidth]{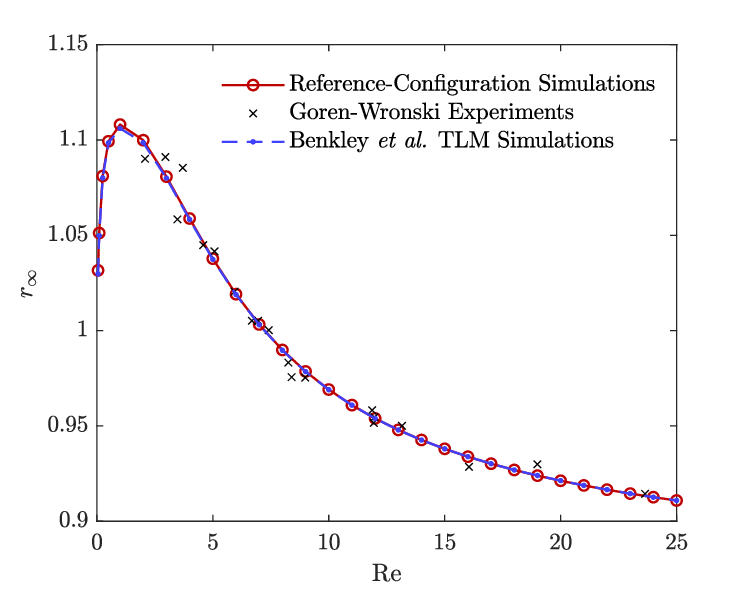}
    \caption{Ratio between terminal radius of jet and radius of the die orifice with respect to $\Rey$ for \textit{Oh} = $\sqrt{3}$. Predictions by simulations with the algorithm proposed in Section \ref{Sec:NumericalMethod} are compared with experimental results found in \cite{goren_shape_1966} and axisymmetric numerical simulations in \cite{benkley_generalization_2026} with the Total Linearization Method (TLM).}
    \label{Fig:DieSwellRatio}
\end{figure}

In Figure \ref{Fig:DieSwellRatio}, the swell ratios predicted by numerical simulation are depicted with respect to the Reynolds number. These data are compared with axisymmetric numerical simulations reported in \cite{benkley_generalization_2026} and with experiments published in \cite{goren_shape_1966}.
Agreement with experiments is observed for a Reynolds number ranging from approximately 2.5 to 25. Predictions are, furthermore, nearly indistinguishable from axisymmetric simulations by the Total Linearization Method for all $\Rey$ shown in Figure \ref{Fig:DieSwellRatio}. This validates the alternative formulation and numerical method presented in Sections \ref{Sec:NavierStokes} and \ref{Sec:NumericalMethod}.

\begin{figure}
    \centering
    \includegraphics[width=0.6\linewidth]{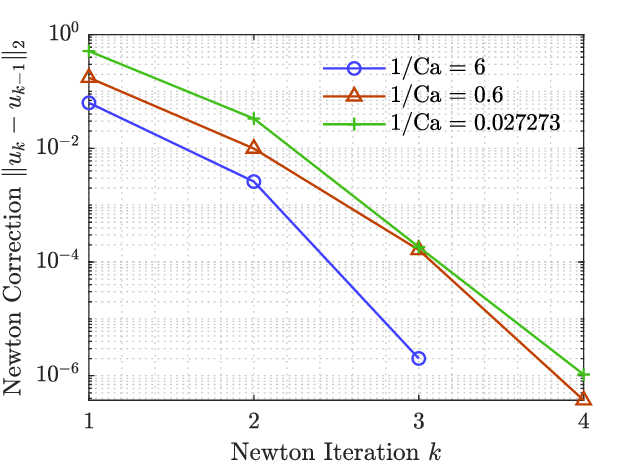}
    \caption{Euclidean norm of Newton correction with respect to Newton iteration of die-swell solve for \textit{Oh} $= \sqrt{3}$. Three different cases of $1/\Ca$ are displayed, each simulated with 8,704,490 degrees of freedom. }
    \label{Fig:DieSwellNewtonConvergence}
\end{figure}
In Figure \ref{Fig:DieSwellNewtonConvergence}, the velocity correction mandated by a Newton update is represented with respect to the Newton iteration. These results show superlinear convergence in three different cases with convergence being achieved in at most four iterations, ascertaining that the inexact Newton method yields a fast-converging nonlinear scheme.

\begin{table}
\centering
    \begin{tabular}{| c | c | c |} 
         $1/\Ca$ & 8,704,490 \# Total DOFs & 68,838,920 \# Total DOFs \\
         \hline
         0.3 & 9.2 mins & 34.7 mins \\
         \hline
         1.2 & 4.9 mins & 20.6 mins \\
         \hline
         6 & 3.5 mins & 18.8 mins \\
         \hline
         100 & 3.8 mins & 15.7 mins \\
         \hline
         10000 & 2.8 mins & 12.2 mins 
    \end{tabular} 
\caption{CPU clock times for die-swell simulations with $Oh = \sqrt{3}$ for two levels of refinement. Computations are parallelized with MPI and distributed across multiple processes of a cluster of computational nodes so as to provide each process with approximately 120,000 degrees of freedom per rank. Each computational node consists of 2 Intel(R) Xeon(R) Platinum 8360Y 2.4 GHz processors with 36 cores each (72 cores per node). Computations with 8,704,490 unknowns are distributed across 72 ranks of one computational node whereas with 68,838,920 unknowns, they are distributed over 576 ranks of 8 total nodes. Clock times include the initial-guess computation.}
\label{Table:DieSwell}
\end{table}
In Table \ref{Table:DieSwell}, the CPU clock times for the entire Newton solution are reported. Calculations are distributed across the computational nodes of a cluster with MPI so as to provide each rank with an average of approximately 120,000 unknowns. For 8,704,490 degrees of freedom, these computations range from 2.8 to 9.2 minutes. For the highest level of refinement, with 68,838,920 unknowns, computations take from 12.2 to 34.7 minutes. 

\begin{figure}
    \centering
    \includegraphics[width=0.6\linewidth]{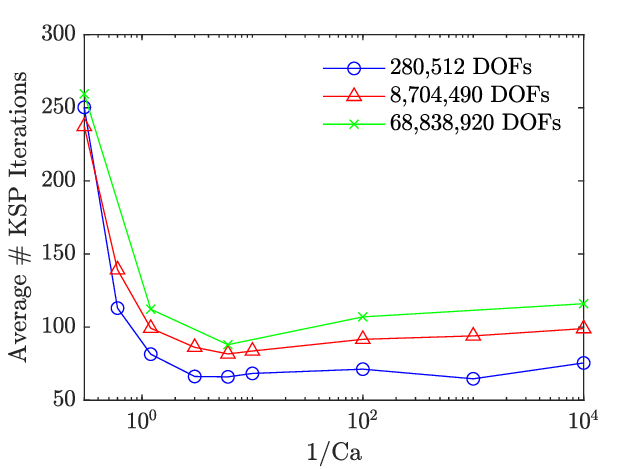}
    \caption{Average number of Krylov-space iterations per Newton iteration of die-swell solve with respect to $1/\Ca$ for \textit{Oh} $= \sqrt{3}$. Results for three levels of refinement are compared. Computations are parallelized and distributed so as to provide each process with approximately $100,000-150,000$ degrees of freedom.}
    \label{Fig:DieSwellKSPvsCa}
\end{figure}

In Figure \ref{Fig:DieSwellKSPvsCa}, the number of Krylov-space iterations averaged over an ensemble of Newton iterations is displayed with respect to $1/\Ca$. The simulations are carried out for three different levels of refinement, leading to linear-system sizes of respectively 280,512, 8,704,490 and 68,838,920 degrees of freedom. Meshes are displayed in Section \ref{Sec:AppendixA}. The iteration count is approximately constant, hovering around 70-110 iterations across nearly four decades of $1/\Ca$ for all three meshes. In comparable cases, this leads to nearly half as many iterations as reported in \cite{benkley_generalization_2026}, addressing precisely the regime where their preconditioner becomes ineffective. These results exhibit effectiveness of the physics-based preconditioner described in Section \ref{Sec:Preconditioning} in the high-surface-tension regime and show little sensitivity to the mesh size. The preconditioner's effectiveness degrades, however, for lower values of $1/\Ca$, where the number of iterations can be observed to diverge as $1/\Ca\searrow 0$. This occurs as a result of the Jacobian limit taken in Equation (\ref{Eq:JacobianLimit}) which breaks down in the low surface-tension regime.

\subsection{Coating Flow} \label{Sec:CoatingFlow}

\begin{figure}
    \centering
    \includegraphics[width=0.85\linewidth]{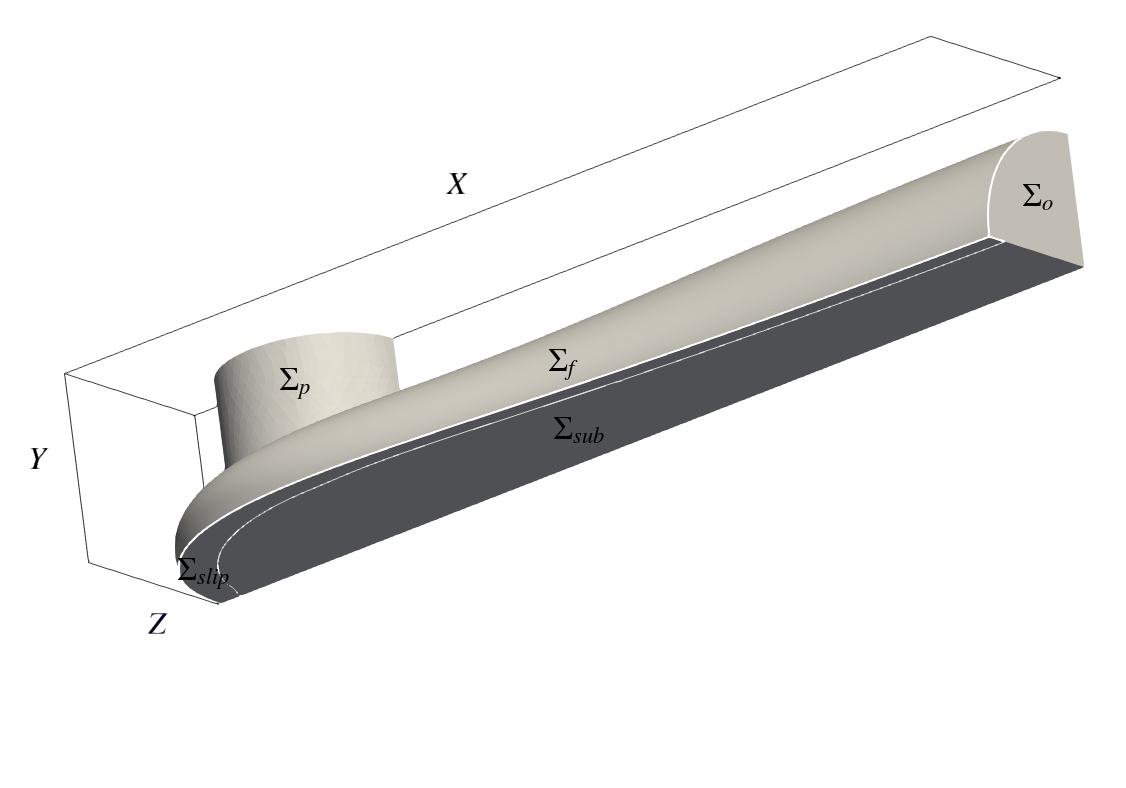}
    \caption{Coating-flow geometry: the liquid enters the numerical domain through the top of the pipe, $\Sigma_p$, before exiting and forming a free interface between liquid and air, $\Sigma_f$. It is deposited onto the solid substrate, which possesses both a no-slip and slip surface, respectively, $\Sigma_\textrm{sub}$ and $\Sigma_\textrm{slip}$. The liquid exits the numerical domain through the outflow surface, $\Sigma_o$. }
    \label{Fig:CoatingFlow}
\end{figure}

In the coating flow considered in Figure \ref{Fig:CoatingFlow}, after entering a numerical inlet with a Poiseuille velocity profile, an incompressible liquid exits a vertical pipe, $\Sigma_p$, downward (along the $y$ axis) into air, thus forming a free interface between liquid and gas. This simulation is inspired by the additive manufacturing process presented in \cite{le_bas_microstructure_2026} and was also modeled previously in \cite{baer_finite_2000}. The free surface, $\Sigma_f$, is subject to the Young-Laplace equation, a shear-free condition and impermeability. The liquid then collides with and coats a solid substrate moving in the horizontal direction (in the $x$ direction) at velocity $U$, satisfying the no-slip condition over the substrate, $\Sigma_\textrm{sub}$. A Navier-slip boundary condition is imposed over a slip region, $\Sigma_\textrm{slip}$, with constant slip length $L_s$, adjacent to the contact line, $\Gamma_{c} = \partial \Sigma_{\textrm{slip}}\cap \partial \Sigma_{\textrm{f}}$, which is subject to an equilibrium contact-angle condition. Further down the substrate, the liquid exits the numerical domain through an outlet, $\Sigma_o$, subject to a homogeneous Neumann boundary condition in the normal direction and zero tangent velocity, imposed as a Dirichlet boundary condition. The steady process is furthermore symmetrical about the $z=0$ plane, $\Sigma_\textrm{sym}$, and is consequently modeled only for $z\leq 0$, with impermeability and zero shear imposed over the symmetry plane.

The characteristic velocity scale is set so as to give a non-dimensionalized flow rate of identically one through the numerical inlet.

Let us introduce the velocity spaces
\[\begin{split}
    V_{D,h} = \{ \v\in [\chi_{2,h}(\Omega_0)]^3 \mid\, &\v|_{\Sigma_{p,0}} = \frac{4}{\pi}(1-x^2-z^2) \begin{pmatrix}0 \\ -1 \\ 0 \end{pmatrix},\v|_{\Sigma_\textrm{sub,0}} = \begin{pmatrix}U\\0\\0 \end{pmatrix},\,v_y|_{\Sigma_\textrm{slip,0}} = 0, \\
    &\,v_z|_{\Sigma_\textrm{sym,0}} = 0, \, v_y|_{\Sigma_{o,0}}= v_z|_{\Sigma_{o,0}} = 0\}
\end{split}\]
and 
\[\begin{split}
    V_{0,h} = \{ \v\in [\chi_{2,h}(\Omega_0)]^3 \mid\, &\v|_{\Sigma_{p,0}} = 0,\v|_{\Sigma_\textrm{sub,0}} = 0,\,v_y|_{\Sigma_\textrm{slip,0}} = 0,\,v_z|_{\Sigma_\textrm{sym,0}} = 0,\,v_y|_{\Sigma_{o,0}}= v_z|_{\Sigma_{o,0}} = 0\}.
\end{split}\]
Moreover, the displacement spaces are defined,
\[\begin{split}
    M = \{ \d\in [H^1(\Omega_0)]^3 \mid \d|_{\Sigma_{f,0}} \in [H^1(\Sigma_{f,0})]^3,\,\d|_{\Sigma_{\textrm{slip},0}} \in [H^1(\Sigma_{\textrm{slip},0})]^3\,\d|_{\Sigma_{\textrm{sub},0}} \in [H^1(\Sigma_{\textrm{sub},0})]^3  \},
\end{split}\]
\[\begin{split}
    M_h = \{ \d\in [\chi_{2,h}(\Omega_0)]^3 \mid\, &\d|_{\Sigma_{p,0}} = 0, \,d_y|_{\Sigma_\textrm{sub,0}} = 0, \,d_y|_{\Sigma_\textrm{slip,0}} = 0 \,d_z|_{\Sigma_\textrm{sym,0}} = 0,\,d_x|_{\Sigma_{o,0}} = 0\},
\end{split}\]
\[\begin{split}
    M_{\textrm{bulk},h} = \{ \d\in M_h \mid \d|_{\Sigma_{f,0}\cup \Sigma_\textrm{slip,0}\cup \Sigma_\textrm{sub,0}} = 0\},
\end{split}\]
\[\begin{split}
    M_{\textrm{sub},h} = \{ \d\in M_h \mid \d|_{\Sigma_{f,0}\cup \Sigma_\textrm{slip,0}} = 0\},
\end{split}\]
\[\begin{split}
    M_{\textrm{slip},h} = \{ \d\in M_h \mid \d|_{\Sigma_{f,0}} = 0\},
\end{split}\]
and 
\[
    M_{f,h} = \{ \d\in M_{h}\mid \d|_{\Gamma_{c,0}} = 0 \}.
\]
In this particular problem, the Laplace-Beltrami Closure Equation (\ref{Eq:Closure}) is tested only over functions $\boldsymbol{\xi}_h$ satisfying $\boldsymbol{\xi}_h|_{\Gamma_{c,0}} = 0$, \textit{i.e.}, 
\[
\int_{\Sigma_{f,0}} \bnabla_S\d_h:\bnabla_S[(I - \frac{\F_0^{-T}\n_0\n_0^T\F_0^{-1}}{|\F_0^{-T}\n_0|^2})\boldsymbol{\xi}_h] dS = 0,\,\forall \xi_h \in M_{f,h}.
\]
This avoids an interplay between boundary terms from the closure equation and the Young-Laplace equation. Equation (\ref{Eq:Closure}) is then complemented with the boundary condition along the contact line,
\[
(\d_h\cdot\t)\t = 0,\textrm{ along $\Gamma_{c}$},
\]
imposed weakly, yielding
\[
\int_{\Gamma_{c}} (\boldsymbol{\xi}_h\cdot \t)(\d_h\cdot\t) dl = 0,\,\forall \boldsymbol{\xi}_h\in M_{h}.
\]
A parametrized and, subsequently, linearized version of this equation is assembled in the Jacobian. 

A Laplace problem is imposed over the slip region,
\begin{equation}
    \begin{cases}- \Delta_S \d_h = 0, \textrm{ over }\Sigma_{\textrm{slip},0} \\
    (\s_r\cdot\bnabla_S)\d_h = 0,\textrm{ along } \partial \Sigma_{\textrm{slip},0} \setminus \Gamma_{c,0}
    \end{cases}
\end{equation}
where the degrees of freedom along $\Gamma_{c,0}$ are determined by the free-surface equations, giving
\[
\int_{\Sigma_{\textrm{slip},0}} \bnabla_S\d_h:\bnabla_S\boldsymbol{\xi}_h \, dS = 0,\,\forall \boldsymbol{\xi}_h\in M_{\textrm{slip},h}.
\]
The homogeneous Neumann boundary condition ensures that the width of the slip region remains approximately constant when deformed. 

A similar Laplace problem is imposed over the substrate, ensuring a constant deformation gradient over the substrate, resulting in
\[
\int_{\Sigma_{\textrm{sub},0}} \bnabla_S\d_h:\bnabla_S\boldsymbol{\xi}_h \, dS = 0,\,\forall \boldsymbol{\xi}_h\in M_{\textrm{sub},h}.
\]
This promotes mesh integrity throughout deformation. The displacement block of the Jacobian thus possesses the block structure,

\[
    J_{dd}=\begin{pmatrix}
        J_{\Omega_0 \Omega_0} & J_{\Omega_0 \Sigma_{\textrm{sub},0}} & J_{\Omega_0 \Sigma_{\textrm{slip},0}}& J_{\Omega_0 \Sigma_{f,0}} \\
        0 & J_{\Sigma_{\textrm{sub},0} \Sigma_{\textrm{sub},0}} & J_{\Sigma_{\textrm{sub},0} \Sigma_{\textrm{slip},0}} & J_{\Sigma_{\textrm{slip},0} \Sigma_{f,0}} \\
        0 & 0 & J_{\Sigma_{\textrm{slip},0} \Sigma_{\textrm{slip},0}} & J_{\Sigma_{\textrm{slip},0} \Sigma_{f,0}} \\
        0 & 0 & 0 & J_{\Sigma_{f,0} \Sigma_{f,0}}
    \end{pmatrix}.
\]
Within a preconditioner application, $J_{dd}$ is approximately inverted via the sequential algorithm described in Section \ref{Sec:DisplacementSolve}. In this particular problem, block $J_{\Sigma_{f,0} \Sigma_{f,0}}$ is inverted via an exact LU solve. Other diagonal blocks are inverted via one multi-grid cycle.

The $M$ norm is made to be
\[\|.\|_M^2 = \|.\|_{H^1(\Omega_0)}^2 + \|.\|_{H^1(\Sigma_{f,0})}^2 + \|.\|_{H^1(\Sigma_{\textrm{slip},0})}^2 + \|.\|_{H^1(\Sigma_{\textrm{sub},0})}^2.\]

Given a reference configuration $\Omega_0$, the problem may be found in weak form below (prior to parametrization).
\begin{equation} \label{Eq:CoatingFlowWeakForm}
    \begin{split}
    &\text{Find $(\u_h,p_h,\lambda_h,\d_h,p_a)\in V_{D,h}\times Q_h \times L_{0,h} \times M_h \times \mathbb{R} $ such that} \\
    &\begin{cases}
        \int_\Omega \{\,\Rey\,\v_h\cdot (\u_h\cdot \bnabla )\u_h  + [\bnabla\u_h+(\bnabla\u_h)^T]:\bnabla\v_h  - p_h(\bnabla\cdot\v_h)\}dV \\
        \text{\textcolor{white}{aaa}} - \int_{\Sigma_f} \lambda_h (\v_h\cdot\n) dS +\frac{1}{L_s} \int_{\Sigma_\textrm{slip}}(\u_h\cdot\v_h) dS =  \frac{1}{L_s} \int_{\Sigma_\textrm{slip}}v_{h,x}\, U dS,\,\hfill\forall \v_h\in V_{0,h} \\
        \int_\Omega -q_h(\bnabla \cdot \u_h)dV = 0,\,\hfill\forall q_h\in Q_h \\
        \int_{\Sigma_f} -\chi_h(\u_h\cdot \n)dS = 0,\,\hfill\forall \chi_h\in L_{0,h} \\
        \int_{\Sigma_{f,0}} \bnabla_S\d_h:\bnabla_S[(I - \frac{\F_0^{-T}\n_0\n_0^T\F_0^{-1}}{|\F_0^{-T}\n_0|^2})\boldsymbol{\xi}_h] dS = 0,\,\hfill\forall \boldsymbol{\xi}_h \in M_{f,h}  \\
        \int_{\Sigma_f} [ (\boldsymbol{\xi}_h\cdot\n)( \Ca \,\lambda_h + p_a) + \bnabla_S\cdot\boldsymbol{\xi}_h  ] dS \\
        \text{\textcolor{white}{aaa}}- \cos{(\theta)} \int_{\Gamma_{c}} (\boldsymbol{\xi}_h\cdot \s_r) dl + \int_{\Gamma_{c}} (\boldsymbol{\xi}_h\cdot \t)(\d_h\cdot\t)dl = 0,\,\hfill\forall \boldsymbol{\xi}_h\in M_{h} \\
        \int_{\Sigma_{\textrm{slip},0}} \bnabla_S\d_h:\bnabla_S\boldsymbol{\xi}_h \, dS = 0,\,\hfill\forall \boldsymbol{\xi}_h\in M_{\textrm{slip},h} \\
        \int_{\Sigma_{\textrm{sub},0}} \bnabla_S\d_h:\bnabla_S\boldsymbol{\xi}_h \, dS = 0,\hfill \forall \boldsymbol{\xi}_h\in M_{\textrm{sub},h} \\
        \int_{\Omega_0} \bnabla\d_h:\bnabla\boldsymbol{\xi}_h\,dV = 0,\,\hfill\forall \boldsymbol{\xi}_h\in M_{\textrm{bulk},h} \\
        U\int_{\Sigma_o}dS = 1
    \end{cases}
    \end{split}
\end{equation}
This problem is parametrized by $\Omega_0$, and then solved with the numerical method explained in Section \ref{Sec:NumericalMethod}. In these simulations, numerical continuation is required to iterate over the spectrum of Capillary numbers investigated here. Newton relaxation is also required in certain cases. Due to sizable deformation of the mesh between certain parameters, remeshing of interior mesh tetrahedra is carried out for highly distorted meshes. 

Nonlinear iteration is considered to have converged when the maximum absolute value of the prescribed Newton correction falls below an absolute tolerance of $10^{-5}$ for $\u_h$ and $\d_h$. Krylov-space iteration is considered to have converged when the Euclidean norm of the residual falls below either a relative tolerance of $10^{-6}$ or an absolute tolerance of $10^{-12}$.

\begin{figure}
    \centering
    \includegraphics[width=0.85\linewidth]{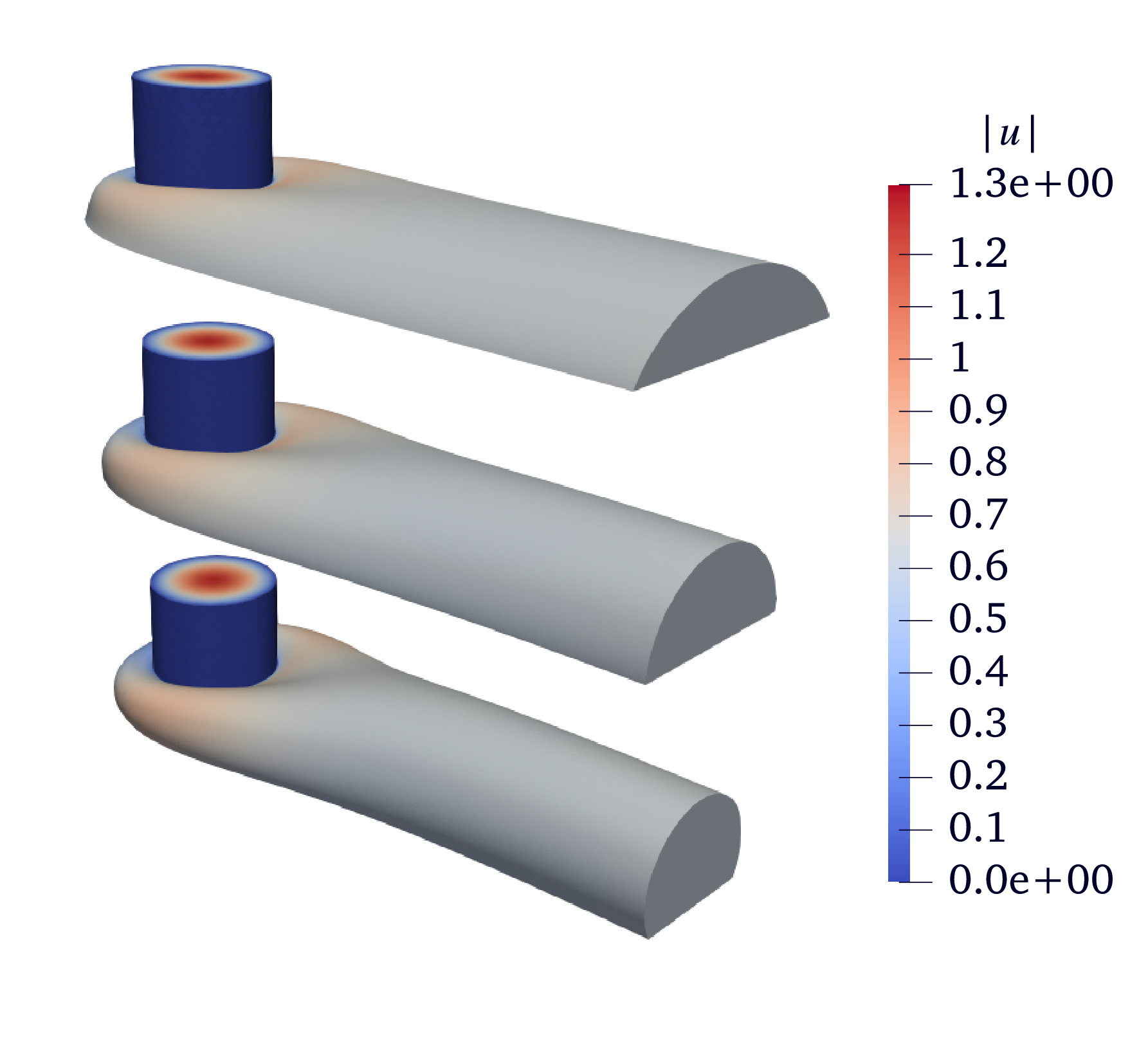}
    \caption{Coating-flow simulation result in the Stokes limit, $\Rey = 0$ for $1/\Ca = 1$, non-dimensionalized substrate velocity $U = 0.6$ and three different contact angles $\theta$. From top to bottom, respectively, $\theta = 60^\circ,90^\circ,120^\circ$. Simulations are viewed from a 3D angle with the velocity magnitude displayed. Simulation renderings are augmented symmetrically about the plane $z=0$.}
    \label{Fig:Diagonalinv_Ca1}
\end{figure}

In Figure \ref{Fig:Diagonalinv_Ca1}, coating-flow results are displayed for simulation parameters $\Rey = 0$, $1/\Ca = 1$, $U = 0.6$ and three different contact angles, respectively, $\theta = 60^\circ$, $90^\circ$ and $120^\circ$. In this case, viscous forces are approximately equal in magnitude to surface-tension forces. Upon exiting the pipe, akin to die-swell flows at low Reynolds numbers (as seen in \cite{benkley_generalization_2026}), the liquid may be observed to bulge out laterally immediately beyond the stick-slip transition due to a similar viscous shearing effect. The amount of substrate wet by the liquid at the numerical outlet, \textit{i.e.} the extent of liquid spreading, is then influenced by the contact angle $\theta$. As the contact angle increases, the energy cost associated with wetting the substrate grows, causing the fluid to wet less of the substrate.

\begin{figure}
    \centering
    \includegraphics[width=0.85\linewidth]{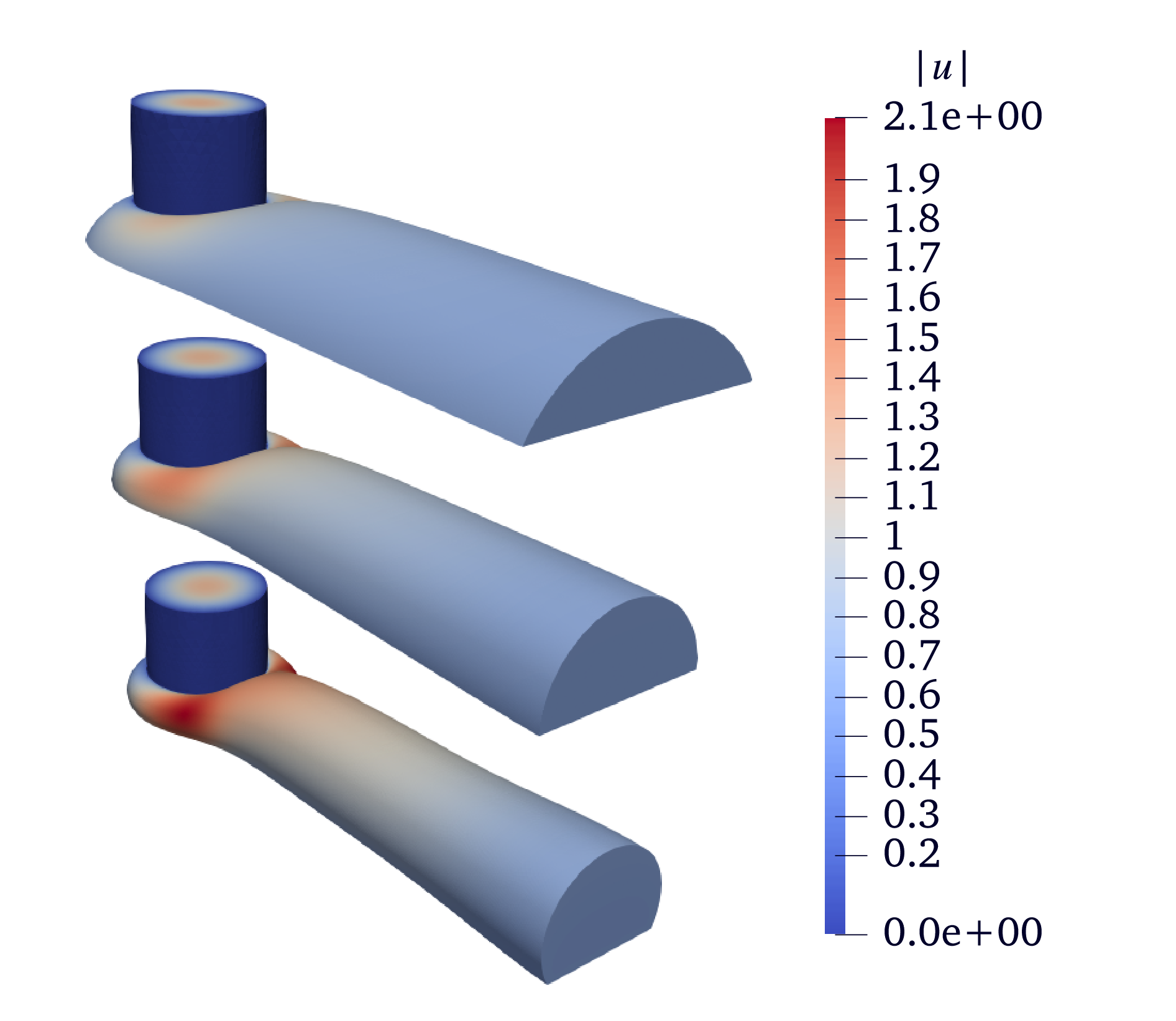}
    \caption{Coating-flow simulation result in the Stokes limit, $\Rey = 0$ for $1/\Ca = 10$, non-dimensionalized substrate velocity $U = 0.6$ and three different contact angles $\theta$. From top to bottom, respectively, $\theta = 60^\circ,90^\circ,120^\circ$. Simulations are viewed from a 3D angle with the velocity magnitude displayed. Simulation renderings are augmented symmetrically about the plane $z=0$.}
    \label{Fig:Diagonalinv_Ca10}
\end{figure}

In Figure \ref{Fig:Diagonalinv_Ca10}, simulation results are depicted for parameters $\Rey = 0$, $1/\Ca = 10$, $U = 0.6$ and again, for $\theta = 60^\circ$, $90^\circ$ and $120^\circ$. The same swelling mechanism as in the case $1/\Ca = 1$ is still perceptible, but is milder in intensity due to an increased relative surface tension. The amount of substrate then wet by the liquid at the numerical outlet is equal to that when $1/\Ca = 1$ due to the following reason. When paired with an area constraint, a stress-free normal boundary condition at the outlet results in a uniform velocity profile while also fixing the pressure gauge to zero at the outlet. This implies in turn that, at the outlet, there are no forces additional to those of capillarity, resulting in the same outlet surface for all Capillary numbers at a fixed substrate velocity.

\begin{figure}
    \centering
    \includegraphics[width=0.85\linewidth]{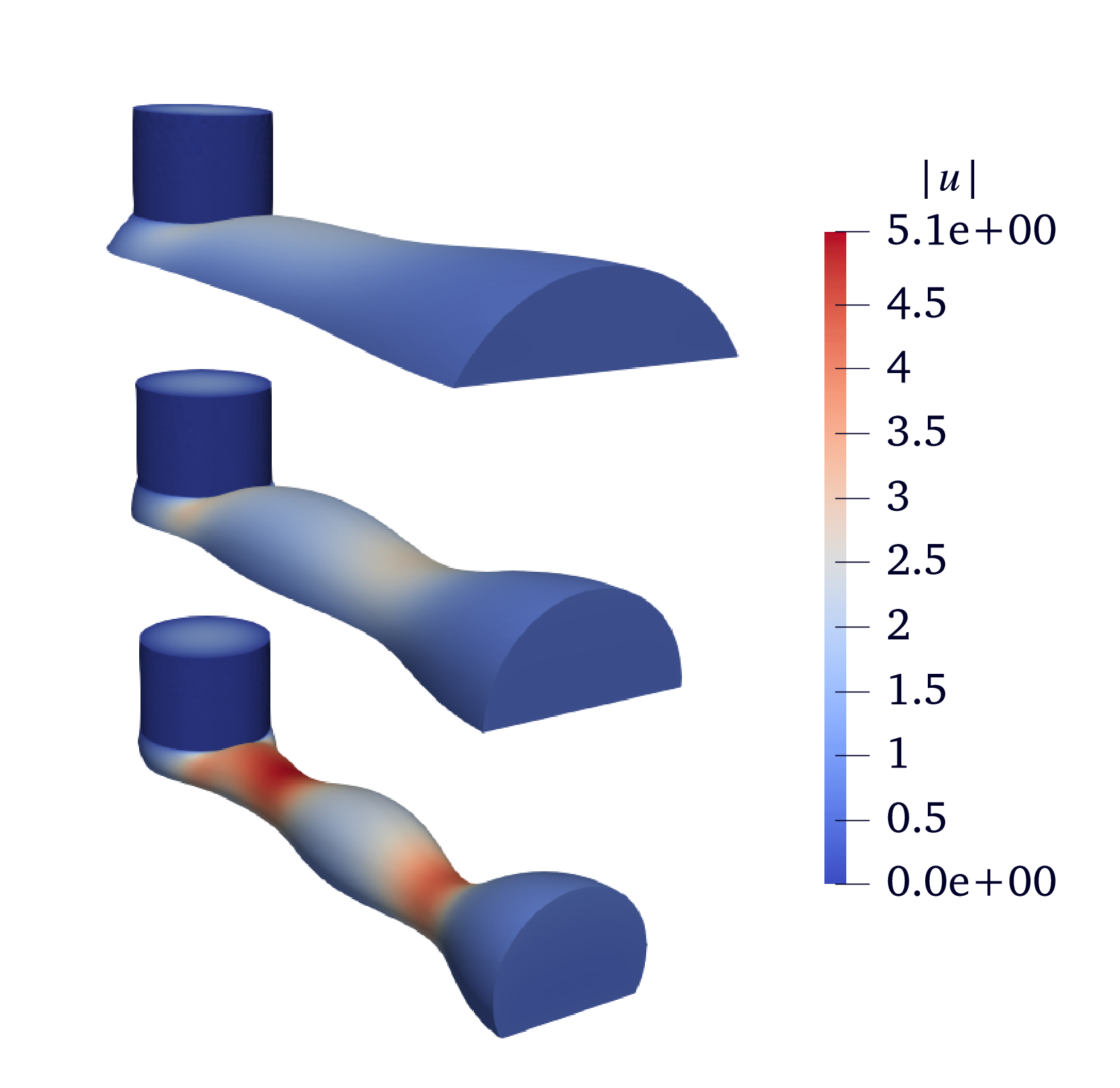}
    \caption{Coating-flow simulation result in the Stokes limit, $\Rey = 0$ for $1/\Ca = 100$, non-dimensionalized substrate velocity $U = 0.6$ and three different contact angles $\theta$. From top to bottom, respectively, $\theta = 60^\circ,90^\circ,120^\circ$. Simulations are viewed from a 3D angle with the velocity magnitude displayed. Simulation renderings are augmented symmetrically about the plane $z=0$.}
    \label{Fig:Diagonalinv_Ca100}
\end{figure}

In Figure \ref{Fig:Diagonalinv_Ca100}, results are shown for simulation parameters $\Rey = 0$, $1/\Ca = 100$, $U = 0.6$ and three different contact angles, respectively, $\theta = 60^\circ$, $90^\circ$ and $120^\circ$. The viscous swelling effect observed in Figure \ref{Fig:Diagonalinv_Ca10} vanishes for $1/\Ca = 100$ due to higher relative surface-tension forces. Indeed, curvature of the leftmost outline of the free surface is fainter than for $1/\Ca = 10$. For $\theta = 90^\circ$ and $\theta = 120^\circ$, conspicuous undulations appear along the longitudinal direction, likely as manifestations of a Rayleigh-Plateau effect described in \cite{schiaffino_formation_1997}. 

\begin{figure}
    \centering
    \includegraphics[width=0.75\linewidth]{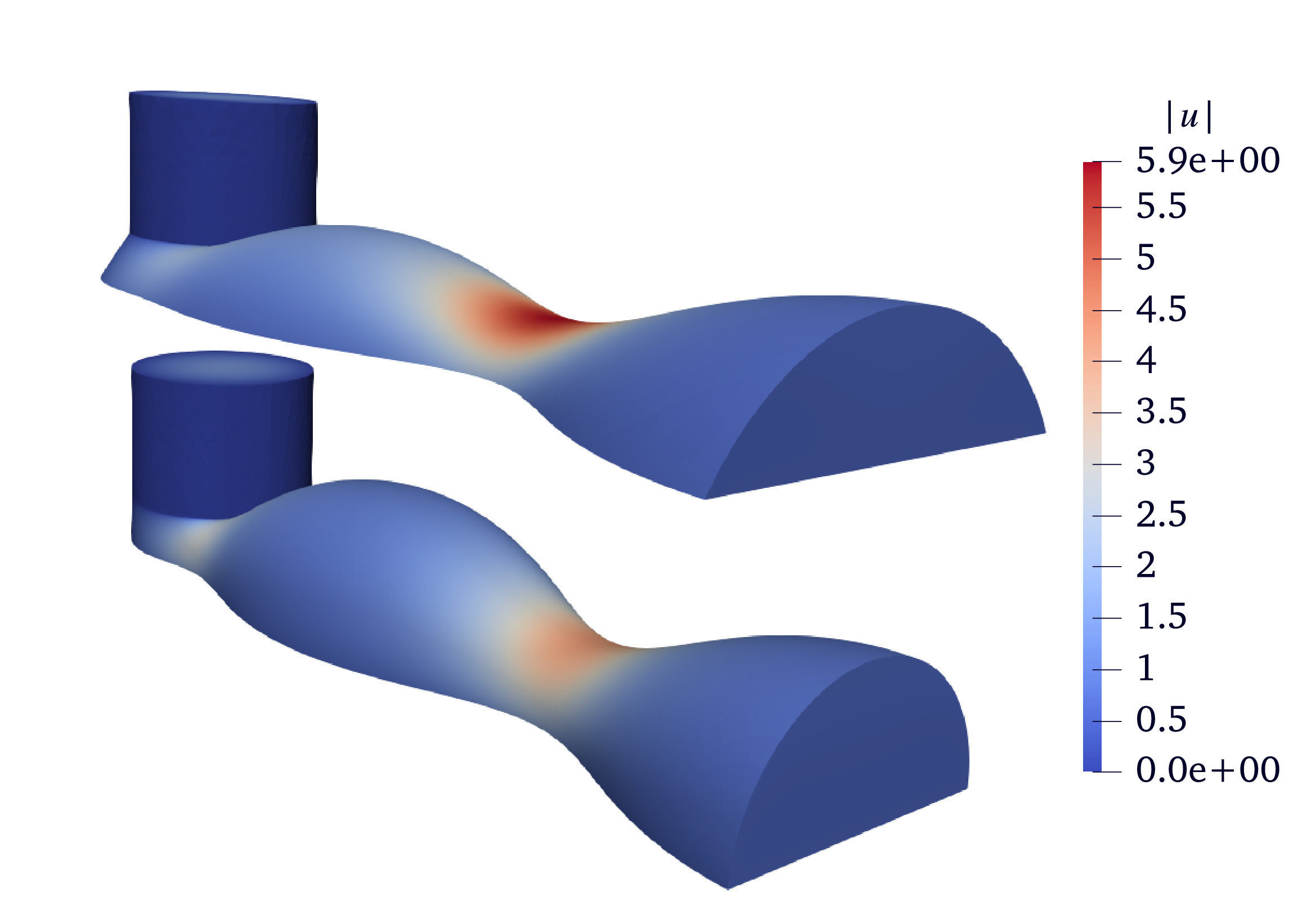}
    \caption{Coating-flow simulation result in the Stokes limit, $\Rey = 0$ for $1/\Ca = 1000$, non-dimensionalized substrate velocity $U = 0.6$ and two different contact angles $\theta$. From top to bottom, respectively, $\theta = 60^\circ,90^\circ$. Simulations are viewed from a 3D angle with the velocity magnitude displayed. Simulation renderings are augmented symmetrically about the plane $z=0$.}
    \label{Fig:Diagonalinv_Ca1000}
\end{figure}
In Figure \ref{Fig:Diagonalinv_Ca1000}, simulation outputs are represented for parameters $\Rey = 0$, $1/\Ca = 1000$, $U = 0.6$ and two different contact angles, respectively, $\theta = 60^\circ$ and $90^\circ$. The contact angle $120^\circ$ is not included as, in this particular case, the elevated surface-tension forces drive the separation of humps into distinct droplets, which the ALE method presented here is not designed to handle. For $\theta = 60^\circ$, prominent bulges of the liquid may be observed, in contrast with cases under milder surface tension. In addition, the bulges for $\theta = 90^\circ$ are further accentuated by the increase in surface tension.

\begin{figure}
    \centering
    \includegraphics[width=0.6\linewidth]{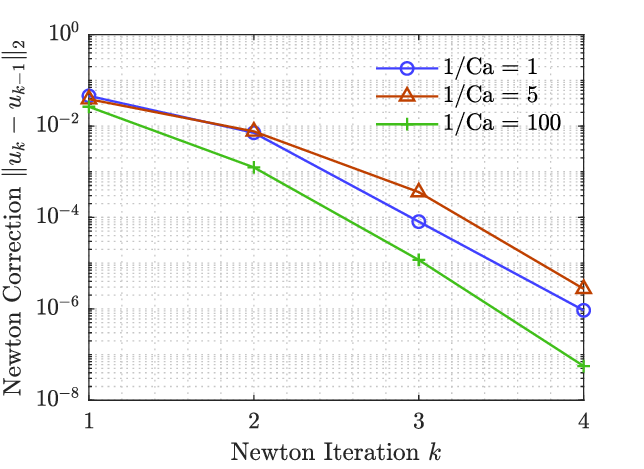}
    \caption{Euclidean norm of Newton correction with respect to Newton iteration of coating-flow solve for $\Rey = 0$, $\theta = 90^\circ$ and $U = 0.6$. Three different cases of $1/\Ca$ are displayed, each simulated with between 2.5 and 6 million degrees of freedom. }
    \label{Fig:CoatingFlowNewtonConvergence}
\end{figure}
In Figure \ref{Fig:CoatingFlowNewtonConvergence}, the velocity correction mandated by a Newton update is represented with respect to the Newton iteration. These results show superlinear convergence in three different cases with convergence being achieved in at most four iterations. Across all simulations, Newton convergence was observed in at most six iterations, ascertaining that the inexact Newton method yields a fast-converging nonlinear scheme.

\begin{figure}
    \centering
    \includegraphics[width=0.6\linewidth]{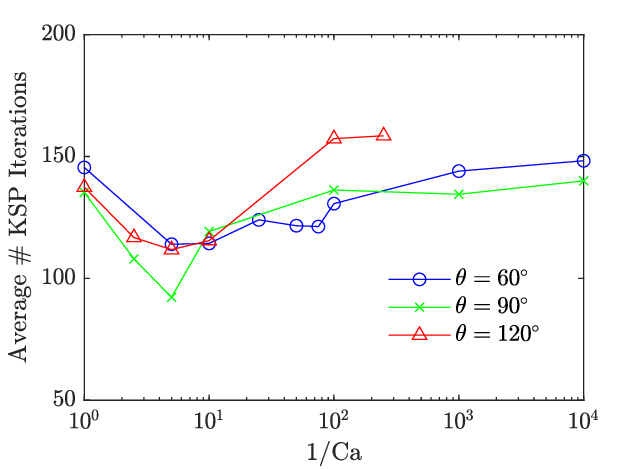}
    \caption{Average number of Krylov-space iterations per Newton iteration of coating-flow solve with respect to $1/\Ca$ for $\Rey = 0$. Between 2.5 and 6 million degrees of freedom are employed in each simulation. Computations are parallelized and distributed across 72 MPI processes.}
    \label{Fig:CoatingFlowKSPvsCa}
\end{figure}

In Figure \ref{Fig:CoatingFlowKSPvsCa}, the number of GMRes iterations averaged over an ensemble of Newton iterations is displayed with respect to $1/\Ca$. The simulations are carried out for meshes leading to between 2.5 and 6 million degrees of freedom. One mesh example is displayed in Section \ref{Sec:AppendixA}. The GMRes iteration count, here, is higher than for the die swell, reaching between approximately 100 and 150 on average. The preconditioner demonstrates adequacy across all four decades of $1/\Ca$ studied here, suggesting suitability specifically in the high-surface-tension regime. As observed in Figures \ref{Fig:Diagonalinv_Ca1}, \ref{Fig:Diagonalinv_Ca10}, \ref{Fig:Diagonalinv_Ca100} and \ref{Fig:Diagonalinv_Ca1000}, these simulations can predict drastically different flow distributions and free surfaces, confirming that preconditioner effectiveness is not unique to one geometry.

\section{Discussion} \label{Sec:Discussion}

The Weak Formulation (\ref{Eq:WeakForm}), an alternative to the classical Weak Form (\ref{Eq:ClassicalWeakForm}), relies on imposing impermeability of the free surface through a Lagrange multiplier $\lambda$. This form is particularly suited to the high surface-tension regime, as surface-tension forces are not present in off-diagonal couplings as in Equation (\ref{Eq:ClassicalWeakForm}) but instead, enter through displacement's diagonal block. To solve the nonlinear weak form, an inexact Newton method is employed, giving rise to a Jacobian matrix that is subsequently factorized into a block preconditioner. The closure equations used for displacement ensure that degrees of freedom over its free surface may be decoupled from those in the bulk within the preconditioner. This allows separate multi-grid hierarchies to be constructed for the different closure equations. The physics-based preconditioner relies on the hypotheses of high surface tension, \textit{i.e.} $\Ca \ll 1$, weak coupling between Lagrange multipliers $p_h$ and $\lambda_h$, and verification of the inf-sup condition by a linear $\lambda_h$ field when paired with the Taylor-Hood elements for $(\u_h,p_h)$. Under these assumptions, with properly tuned multi-grid approximations of sparse inverses, the preconditioner derived in Section \ref{Sec:Preconditioning} is expected to lead to Krylov-space solves characterized by moderate iteration counts for high surface tension and by low sensitivity to the mesh size. 

Simulation results for the die-swell problem result in swell-ratio predictions that match those made in a previous study in \cite{benkley_generalization_2026}, verifying the alternative Weak Form (\ref{Eq:WeakForm}). The Newton solves also all terminated in at most four iterations, confirming that the Jacobian inexactness does not severely degrade the convergence rate of Newton's algorithm. As shown in Table \ref{Table:DieSwell}, the total clock times are at most 35 minutes for 68.8 million unknowns, and 9.2 minutes for 8.7 million unknowns. This makes this method competitive with existing simulation techniques.

The present die-swell simulations yield a lower number of iteration counts than in \cite{benkley_generalization_2026}, in the high surface-tension regime. This iteration count lies between 70 to 110 over four orders of magnitude of $1/\Ca$, from $1$ to $10000$. Only around $\Ca \sim 2$, well after violating the assumption $\Ca\ll 1$, does performance of the preconditioner sizably deteriorate, in contrast with the preconditioner presented in \cite{benkley_generalization_2026}, which performs well in this regime. This occurs due to the neglected coupling in Equation (\ref{Eq:JacobianLimit}) between $\lambda_h$ and $\d_h$. These results verify that the block structure of the preconditioner is indeed suited to high surface-tension flows. The iteration counts are also only mildly sensitive to the mesh size. Indeed, in the worst case, the number of iterations increases by a factor of 1.5 between simulations with 280,512 and 68,838,920 degrees of freedom, \textit{i.e.} an increase by a factor of 245 in the number of unknowns. This suggests both that Schur complement approximations in Section \ref{Sec:PreconditionerApproximation} verify spectral equivalence and that $p_h$ and $\lambda_h$ are indeed weakly coupled. 

The coating-flow simulation analyzes flows more complex than in the die-swell case. Indeed, in Figures \ref{Fig:Diagonalinv_Ca1}, \ref{Fig:Diagonalinv_Ca10}, \ref{Fig:Diagonalinv_Ca100} and \ref{Fig:Diagonalinv_Ca1000}, widely different free surfaces emerge from simulation with, for example, significant swelling at the pipe exit for relatively low surface tension or, for high surface tension, the development of undulations which intensify as $1/\Ca$ grows. Undulations were not observed in the earlier simulations in \cite{baer_finite_2000}, as the largest value of $1/\Ca$ investigated was 2. We believe that these undulations occur as a result of a Rayleigh-Plateau mechanism detailed in \cite{schiaffino_formation_1997, bostwick_static_2018}.

In \cite{schiaffino_formation_1997,bostwick_static_2018}, an infinite, partial cylinder with free contact lines was shown to be linearly unstable to longitudinal perturbations for all contact angles. This occurs due to the fact that, for free contact lines, an increase in cross-sectional area is accompanied by a decrease in mean pressure whereas the opposite happens for shrinking cross sections. This leads to a pump mechanism where smaller cross sections nourish the larger ones. In experiments by Schiaffino and Sonin with free contact lines \cite{schiaffino_formation_1997}, liquid humps ultimately split up into respective sessile droplets. We believe that the ripples observed here arise as a result of a similar, saturated instability. In practice, however, our steady simulations are limited by the domain length and the global constraint which sets the outlet area. This has a direct consequence on the numerical solution, since segregated droplets do not necessarily meet the outlet-area condition in experiments, which may explain why our simulations do not always result in separated droplets for high surface tension. As a side note, in additive manufacturing techniques such as in \cite{le_bas_microstructure_2026}, these Rayleigh-Plateau instabilities are not always observed, as they are often suppressed by solidification, which is not modeled here. 

In Figure \ref{Fig:CoatingFlowKSPvsCa}, GMRes iteration counts for the coating flow reveal a relatively stable number of iterations across the investigated interval of $1/\Ca$, $[1,10000]$, with peak performance achieved near $1/\Ca = 10$ for all three contact angles, $60^\circ$, $90^\circ$ and $120^\circ$. For larger values of surface tension, Figures \ref{Fig:Diagonalinv_Ca100} and \ref{Fig:Diagonalinv_Ca1000} exhibit more complex flows due to undulations which may lead to an increase in the number of iterations. Nevertheless, the increase at high $1/\Ca$ seems to be plateauing as $1/\Ca \nearrow$. The iteration counts also grow when $1/\Ca \searrow 1$, as is the case in the die-swell problem due to the neglected $\d_h-\lambda_h$ coupling. 

The number of GMRes iterations is larger for the coating flow than for the die-swell problem because of the increased flow complexity. Indeed, in the die-swell problem, in addition to being axisymmetric, the free surface remains relatively close to its initial guess of a cylinder. The longitudinal component of displacement, $d_z$, is furthermore made to be zero. These facts simplify the numerical solution, limiting variation of scalar values $p_h$ and $\lambda_h$ along the azimuthal direction. In addition, the coating flows possess free contact lines and significant changes in the free surface, in particular for $1/\Ca \sim 1$, where there are both concave and convex free-surface regions. This added flow complexity can lead to higher numbers of iterations.

There are also subtle fluctuations in the iteration counts in the coating-flow case which do not appear for the die swell. This occurs as a result of mesh changes between different simulations. It was found that the number of GMRes iterations is sensitive to the quality of the mesh. As a result, throughout the numerical-continuation trajectories undergone over $1/\Ca \in [1,10000]$, the interior elements of the mesh were regenerated by Gmsh \cite{geuzaine_gmsh_2009} as necessary, using the mesh boundary facets as input. Prior to regeneration, boundary facets were also optimized with Mmg \cite{dapogny_three-dimensional_2014} to mitigate mesh degradation due to deformation. This leads to potential differences in mesh quality across coating-flow simulations which can cause fluctuations in iteration counts.

In Figure \ref{Fig:CoatingFlowNewtonConvergence}, the evolution of the Newton correction over the nonlinear iteration confirms that, in coating-flow simulations, convergence remains fast despite inexactness of the Jacobian. Convergence is observed in at most six iterations.

\section{Conclusion}

An efficient, ALE method is proposed for the solution of steady 3D free-surface flows in the high surface-tension regime, demonstrating how the infinite surface-tension limit can be leveraged to produce robust and efficient preconditioners.

The method is predicated on an alternative weak formulation of free-surface flows that imposes impermeability of the free surface via the Lagrange multiplier method. Paired with an inexact Newton method, this ensures that surface-tension forces manifest themselves through the diagonal displacement block rather than through the off-diagonal $\u_h-\d_h$ block. This bypasses the need to approximate capillarity in dense Schur complements when using a block factorization, thereby improving the treatment of capillary contributions within block preconditioners. Closure equations are chosen such that, within the preconditioner, displacement degrees of freedom are decoupled from those in the bulk. This allows the treatment of these subsets of degrees of freedom by separate multi-grid hierarchies.

The resulting physics-based preconditioner exploits the structure of the Jacobian in the high-surface-tension limit and incorporates a fractional Sobolev approximation of the impermeability Schur complement. Under the assumed inf-sup stability of the associated Lagrange multiplier, the spectral-equivalence argument provides a theoretical justification for this approximation. This illustrates how the functional setting of a constraint can be exploited to construct preconditioners for the corresponding Schur complement, an approach that may also be relevant to other constrained free-boundary problems.

Numerical simulations of two 3D free-surface flows demonstrate effectiveness of the block preconditioning strategy over four orders of magnitude of the Capillary number, with GMRes iteration counts ranging between 75 and 150 for the investigated configurations. Moreover, the convergence of the inexact Newton scheme indicates that the approximations introduced in the linearization do not severely compromise nonlinear convergence for the cases considered. Parallel computations involving up to 68 million unknowns, completed within runtimes on the order of tens of minutes, demonstrate the feasibility of large-scale three-dimensional simulations in the surface-tension-dominated regime.

While these results establish the practical effectiveness of the proposed method for the steady configurations considered, several extensions remain of interest. In particular, extending the preconditioning strategy to time-dependent free-surface flows would require accounting for the additional couplings introduced by temporal discretization. Further improvements in parallel scalability, particularly in the multi-grid solution of the individual blocks, would also be desirable for simulations involving substantially larger numbers of degrees of freedom.

\section*{Acknowledgments}

This research project is sponsored by an Advanced Grant from the Swiss National Science Foundation, Grant No. TMAG-2\_209328. We would like to acknowledge Professor Andreas Mortensen for helpful comments on the manuscript, Philipp Weder for fruitful discussions and blackboard sessions, and the SCITAS platform of EPFL for computational resources.

\appendix

\counterwithin{figure}{section}

\section{Appendix} 

\subsection{Die-Swell Meshes} \label{Sec:AppendixA}

\begin{figure}
    \centering
    \subfigure[]{\includegraphics[width=0.32\linewidth]{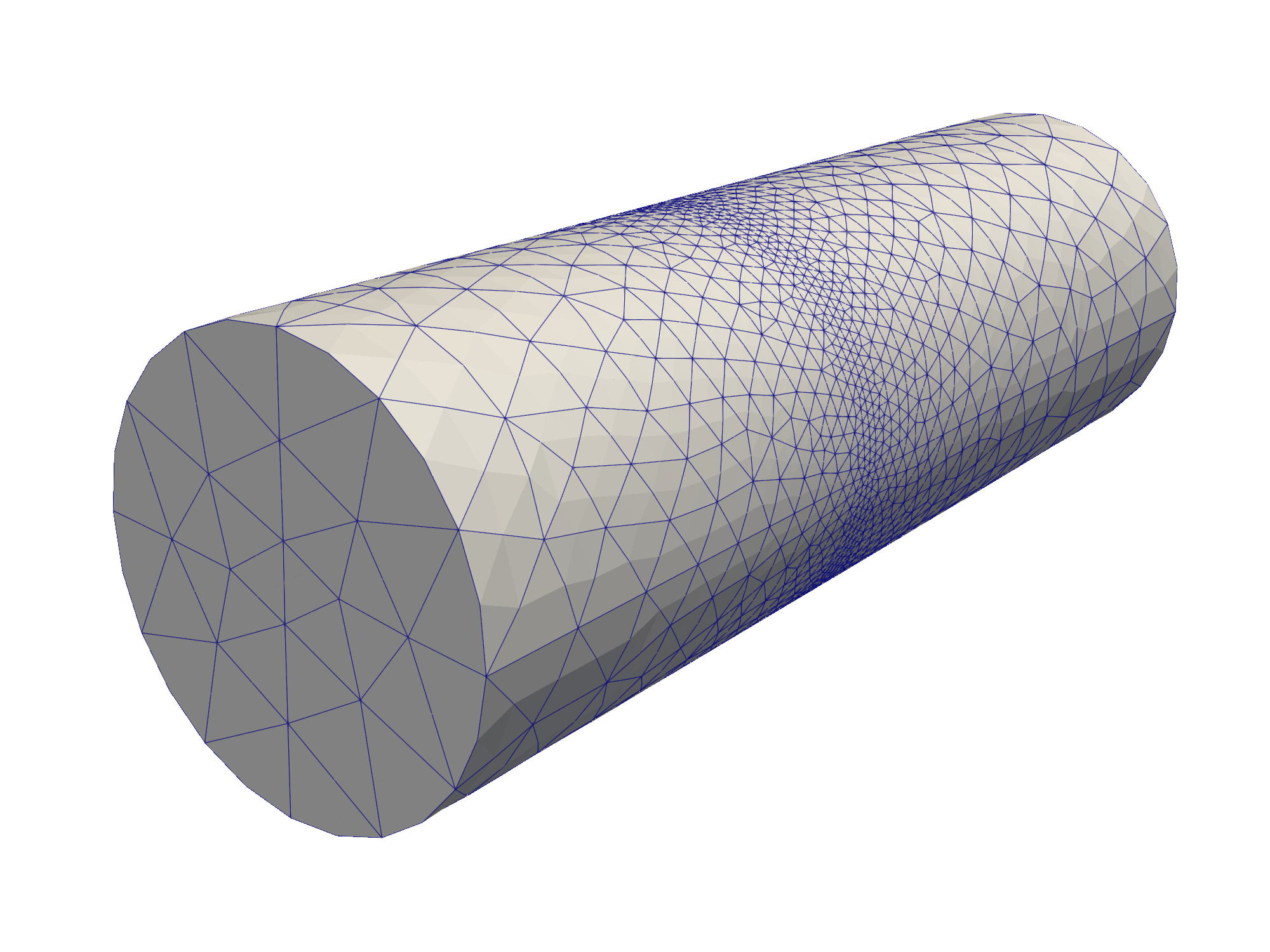}}
    \subfigure[]{\includegraphics[width=0.32\linewidth]{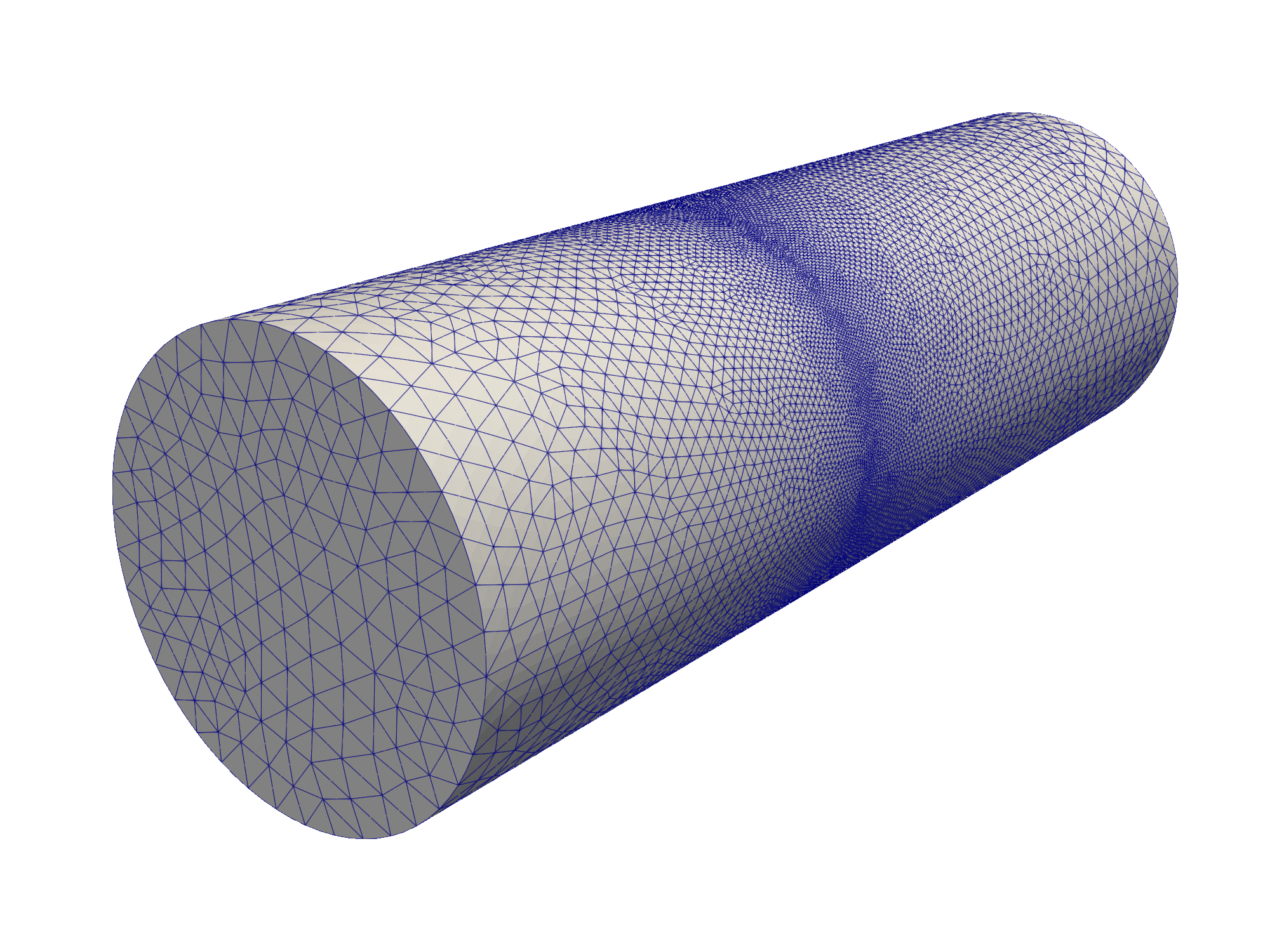}}
    \subfigure[]{\includegraphics[width=0.32\linewidth]{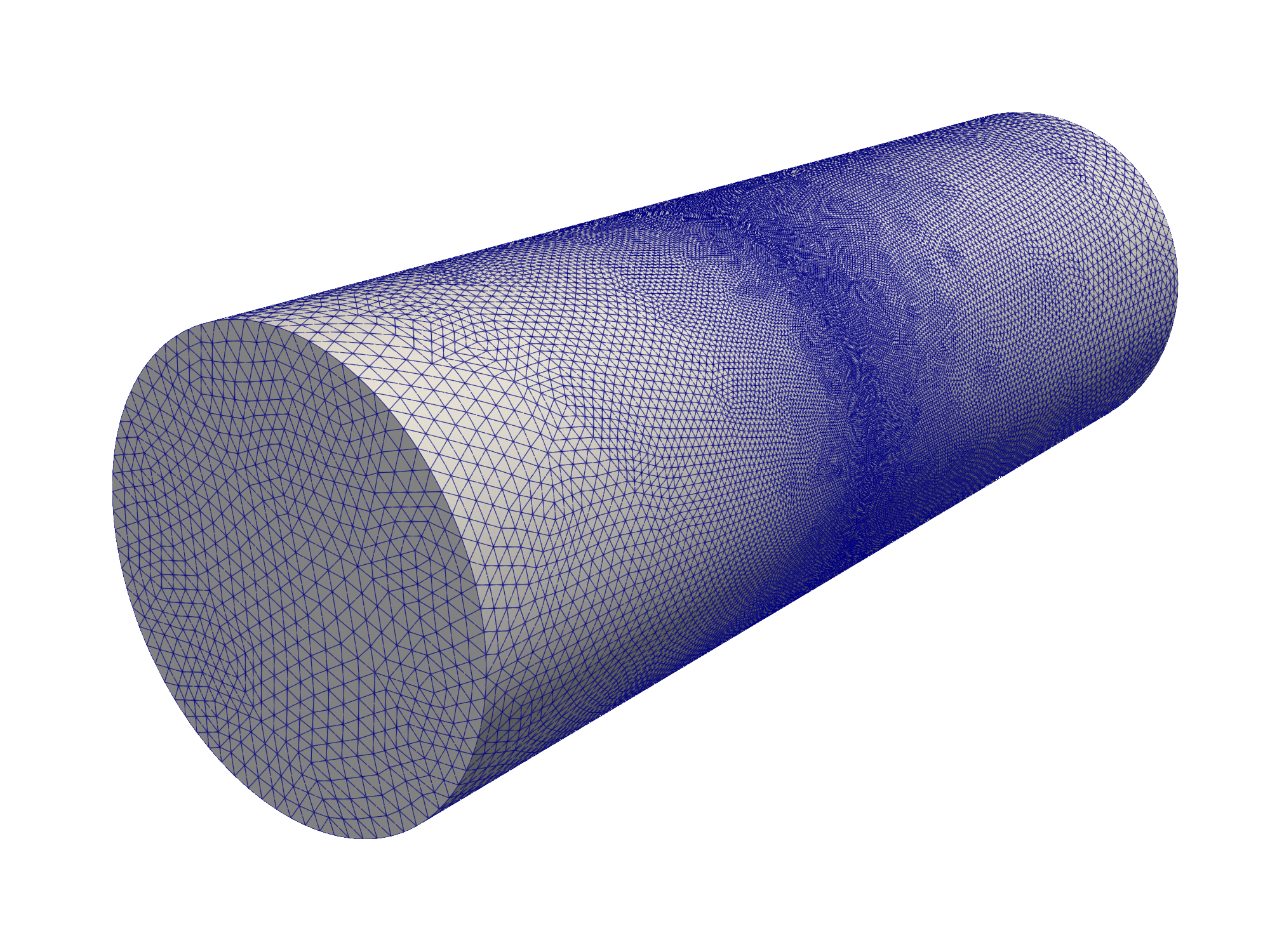}}
    \caption{Spatial discretization for a 3D Cartesian die-swell simulation that led to 280,512, 8,704,490 and 68,838,920 degrees of freedom in (a), (b) and (c), respectively. }
    \label{fig:DieSwellMeshes}
\end{figure}

In Section \ref{Sec:DieSwell}, simulations are carried out with three different levels of refinement, the coarsest mesh leading to 280,512 total degrees of freedom, the medium mesh to 8,704,490 and the finest to 68,838,920. These three meshes are displayed in Figure \ref{fig:DieSwellMeshes}.

\subsection{Coating-Flow Meshes} \label{Sec:AppendixB}

\begin{figure}
    \centering
    \subfigure[]{\includegraphics[width=0.49\linewidth]{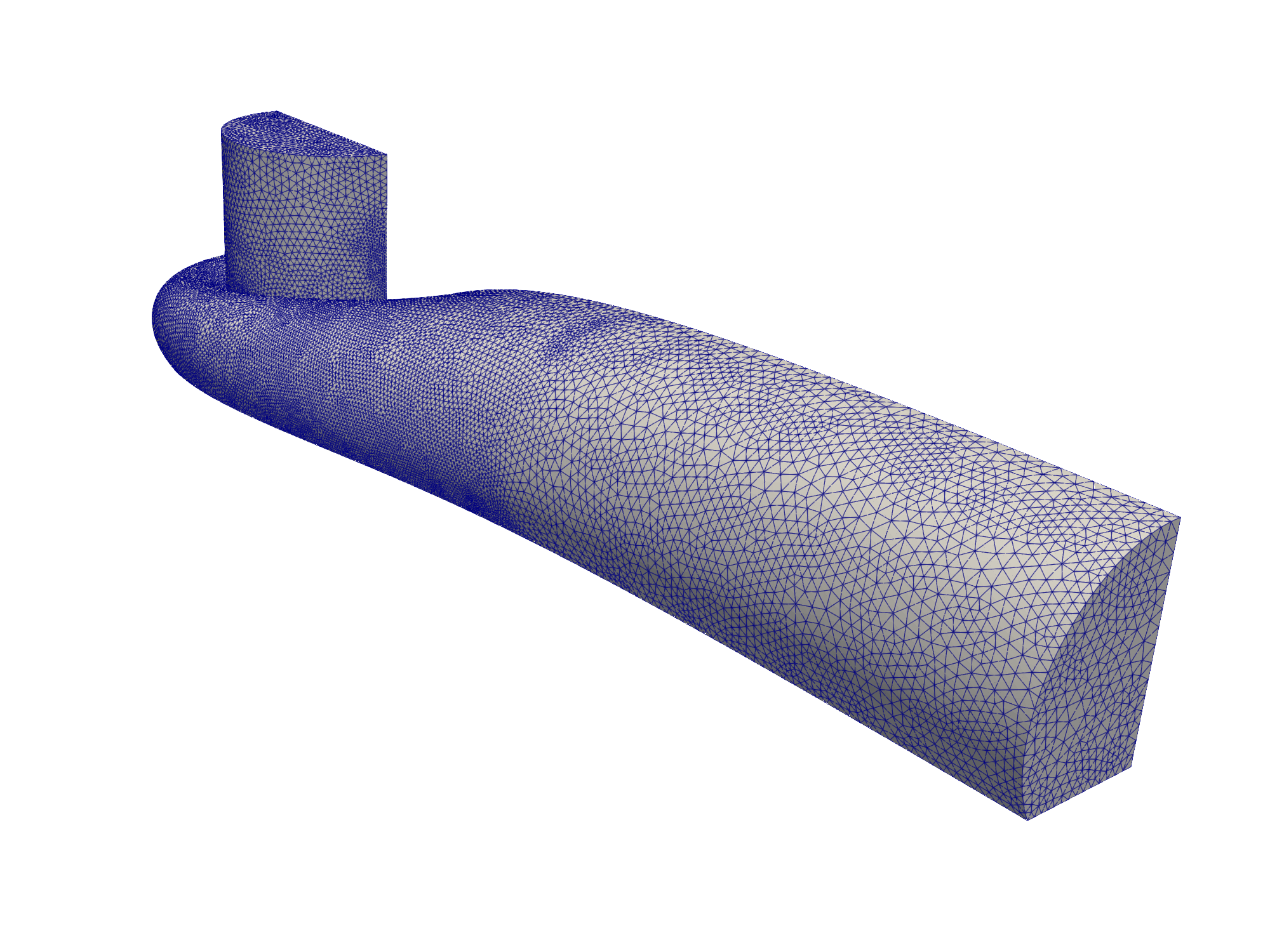}}
    \subfigure[]{\includegraphics[width=0.49\linewidth]{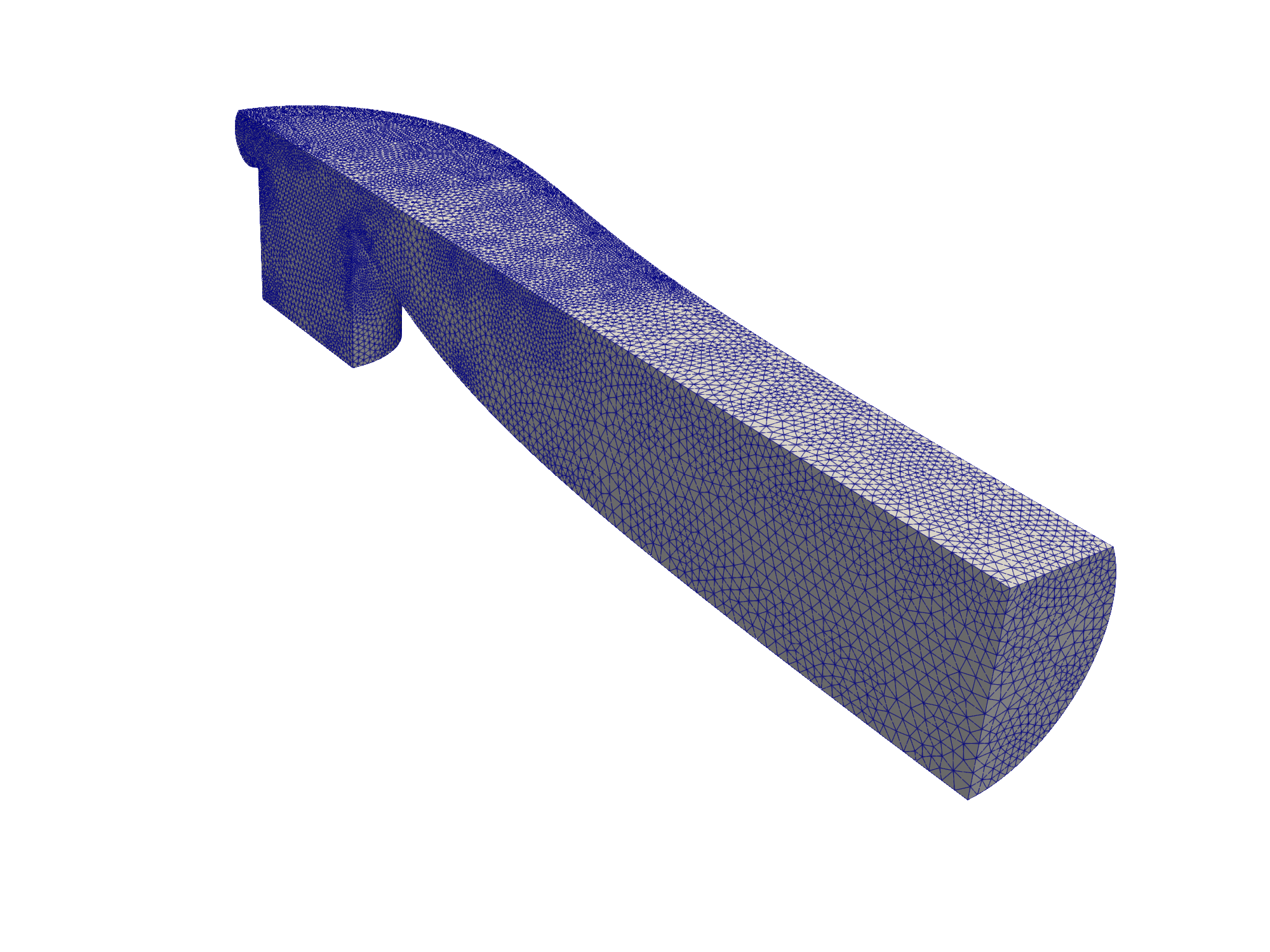}}
    \caption{Spatial discretization for a coating-flow simulation that led to 3,566,021 degrees of freedom, as viewed from the top in (a) and bottom in (b).}
    \label{fig:CoatingFlowMesh}
\end{figure}

In Section \ref{Sec:CoatingFlow}, simulations are carried out for between 2.5 and 6 million degrees of freedom. An example of one mesh is given in Figure \ref{fig:CoatingFlowMesh}.

\newpage

\printbibliography

\end{document}